# A DUALITY APPROACH FOR THE WEAK APPROXIMATION OF STOCHASTIC DIFFERENTIAL EQUATIONS

By Emmanuelle Clément, Arturo Kohatsu-Higa[1]
and Damien Lamberton

*Université de Marne-la-Vallée, Osaka University
and Université de Marne-la-Vallée*

In this article we develop a new methodology to prove weak approximation results for general stochastic differential equations. Instead of using a partial differential equation approach as is usually done for diffusions, the approach considered here uses the properties of the linear equation satisfied by the error process. This methodology seems to apply to a large class of processes and we present as an example the weak approximation of stochastic delay equations.

**1. Introduction.** The Euler scheme for stochastic differential equations is widely used in applications as it is easy to compute. The Euler scheme can be easily generalized to a variety of stochastic equations beyond the framework of diffusion equations, in particular Volterra SDEs, delay SDEs, anticipating SDEs and nonlinear SDEs.

On the other hand, the theoretical properties of the Euler scheme are mostly studied for the diffusion case as most of the results available so far are in this framework. In some cases, extensions to other similar equations are straightforward but in other cases, additional nontrivial work is required. For example, see [8] for extensions to semimartingales, and [1, 11] for approximations of an irregular functional of a diffusion which is approached using a Euler type scheme. It is also well known that the definition of an extension of the Euler scheme for delay type systems is straightforward but the technical results on the weak rate of convergence are limited. See [4, 6, 9].

In this article we propose a generalization of the theory of weak approximations which studies the rate of convergence of the Euler scheme considered

---

Received December 2004; revised October 2005.
[1]Supported in part by Grants BFM 2003-03324 and BFM 2003-04294.

in law. This generalization finds as an application the weak rate of convergence of smooth functionals of general delay type systems and also covers, with a further study of the Malliavin covariance matrix, the case of irregular functions of the solution of the stochastic equation.

The main idea is to change completely the approach used until now to prove weak approximation rate results. This new idea, which uses the whole path of the process under study rather than the partial differential equation associated to the problem, should allow to obtain various other straightforward generalizations of results of the weak rate of convergence.

In order to describe our approach roughly, let $X$ denote the solution of a stochastic equation and $\bar{X}$ the Euler scheme associated to it. The problem of weak rate of convergence consists in finding the rate at which $E(f(X) - f(\bar{X}))$ converges to zero for various classes of functionals $f$. The optimal rate is the step size of the scheme even though the equations considered may differ.

The classical proof of this result for diffusions is based on the associated partial differential equation, that is, $Ef(X)$ has through the Feynman–Kac formula an interpretation using PDEs. This is the important point in the classical approach which is not used in our approach. In the case of some stochastic equations, if $f$ is regular enough, the proof is similar if the associated PDE exists. If $f$ is an irregular function, then the issue of the nondegeneracy of the Malliavin covariance matrix of the Euler scheme becomes an important issue as has been shown in [1, 11], but this extension is nontrivial.

In this article we propose a completely different method to prove weak approximation results based on a pathwise approach. That is, we use the mean value theorem to rewrite $f(X) - f(\bar{X}) = \int_0^1 f'(aX + (1-a)\bar{X}) \, da(X - \bar{X})$. Then, we derive a linear equation satisfied by $Y = X - \bar{X}$. When this equation can be explicitly solved, which seems to be true only for diffusions, one can obtain the rate of convergence by using the duality property of stochastic integrals. This methodology was first introduced by Kohatsu-Higa and Pettersson [7] and used in Gobet and Munos [5]. It seems to be quite general except for the explicit expression for $Y$ which can be done only in the case of diffusions. This article presents a general framework to analyze weak approximations in stochastic equations. In particular we solve the problem without having an explicit expression for the solution of stochastic linear equations, by using a duality argument. This duality formula (see Section 3) shows explicitly the weak error as a by-product of the expectation of an error process (called $G$ in Section 3). To finish the proof one has to use the duality formula for stochastic integrals. Therefore our approach works mostly for stochastic equations with regular coefficients. For this reason we have to study the stochastic derivatives of the solution process. It should be emphasized that this approach applies to regular functionals of the process $X$ and not only to functions of the value of $X$ at a fixed time $t$.

Furthermore, the framework introduced here also extends naturally to the case of irregular functions $f$. That is, one uses the integration by parts formula of Malliavin calculus to regularize the function $f$. We believe that this approach for the irregular case follows naturally from the regular case. It also explains clearly that to obtain the result for the irregular case is just a matter of studying the nondegeneracy of the limiting stochastic process and not the approximating process.

In order to show what are the elements in each theorem we also consider as an example a delay type equation which does not have a clear associated PDE in order to extend the classical proof. The general framework is directly applied, which gives the weak rate of convergence.

This article is organized as follows. In Section 2 we present the example of a diffusion process to introduce the methodology. Section 3 contains the main results in the general framework and an application to the weak rate of convergence of approximations of delay equations for regular functions. These results are then extended to irregular functions in Section 4.

**2. The case of one-dimensional diffusions.** To clarify the methodology, we consider a smooth function $\sigma$ and a Wiener process $W$ and a real diffusion process

$$X_t = x + \int_0^t \sigma(X_s)\,dW_s, \qquad t \in [0, T].$$

Its Euler approximation is given by

$$\bar{X}_t = x + \int_0^t \sigma(\bar{X}_{\eta(s)})\,dW_s, \qquad t \in [0, T],$$

where $\eta(s) = kT/n$ for $kT/n \le s < (k+1)T/n$. The error process $Y = X - \bar{X}$ satisfies

$$Y_t = \int_0^t (\sigma(X_s) - \sigma(\bar{X}_{\eta(s)}))\,dW_s$$

$$= \int_0^t \int_0^1 \sigma'(aX_s + (1-a)\bar{X}_{\eta(s)})\,da\,(X_s - \bar{X}_{\eta(s)})\,dW_s;$$

this can be written as

(1) $$Y_t = \int_0^t \sigma_1(s)Y_s\,dW_s + G_t, \qquad 0 \le t \le T,$$

with

$$\sigma_1(s) = \int_0^1 \sigma'(aX_s + (1-a)\bar{X}_{\eta(s)})\,da,$$

$$G_t = \int_0^t \sigma_1(s)(\bar{X}_s - \bar{X}_{\eta(s)})\,dW_s = \int_0^t \sigma_1(s)\sigma(\bar{X}_{\eta(s)})(W_s - W_{\eta(s)})\,dW_s.$$

In this simple case we have an explicit expression for $Y_t$:

$$Y_t = \mathcal{E}_t \int_0^t \mathcal{E}_s^{-1}(dG_s - \sigma_1(s)\,d\langle G, W\rangle_s),$$

where $\mathcal{E}$ is the unique solution of

$$\mathcal{E}_t = 1 + \int_0^t \sigma_1(s)\mathcal{E}_s\,dW_s.$$

Finally we obtain

$$Y_t = \mathcal{E}_t \int_0^t \mathcal{E}_s^{-1}\sigma_1(s)\sigma(\bar{X}_{\eta(s)})(W_s - W_{\eta(s)})\,dW_s$$

$$- \mathcal{E}_t \int_0^t \mathcal{E}_s^{-1}\sigma_1(s)^2\sigma(\bar{X}_{\eta(s)})(W_s - W_{\eta(s)})\,ds.$$

Now let $f$ be a regular function. We are interested in the determination of the rate of convergence of $\mathbb{E}f(X_T)$ to $\mathbb{E}f(\bar{X}_T)$. We first write the difference

$$\mathbb{E}f(X_T) - \mathbb{E}f(\bar{X}_T) = \mathbb{E}\int_0^1 f'(aX_T + (1-a)\bar{X}_T)\,da\,Y_T.$$

Replacing $Y_T$ by its expression, we obtain with the additional notation $F^h = \int_0^1 f'(aX_T + (1-a)\bar{X}_T)\,da$,

$$\mathbb{E}F^hY_T = \mathbb{E}F^h\mathcal{E}_T \int_0^T \mathcal{E}_s^{-1}\sigma_1(s)\sigma(\bar{X}_{\eta(s)})(W_s - W_{\eta(s)})\,dW_s$$

(2)

$$- \mathbb{E}F^h\mathcal{E}_T \int_0^T \mathcal{E}_s^{-1}\sigma_1(s)^2\sigma(\bar{X}_{\eta(s)})(W_s - W_{\eta(s)})\,ds.$$

Applying the duality for stochastic integrals, this gives

$$\mathbb{E}F^hY_T = \mathbb{E}\int_0^T D_s(F^h\mathcal{E}_T)\mathcal{E}_s^{-1}\sigma_1(s)\sigma(\bar{X}_{\eta(s)})(W_s - W_{\eta(s)})\,ds$$

$$- \mathbb{E}F^h\mathcal{E}_T \int_0^T \mathcal{E}_s^{-1}\sigma_1(s)^2\sigma(\bar{X}_{\eta(s)})(W_s - W_{\eta(s)})\,ds,$$

where $D$ denotes the stochastic derivative. Consequently, the difference $\mathbb{E}f(X_T) - \mathbb{E}f(\bar{X}_T)$ can be written as

$$\mathbb{E}f(X_T) - \mathbb{E}f(\bar{X}_T) = \mathbb{E}\int_0^T U_s^h(W_s - W_{\eta(s)})\,ds,$$

with

$$U_s^h = (D_s(F^h\mathcal{E}_T) - F^h\mathcal{E}_T\sigma_1(s))(\mathcal{E}_s^{-1}\sigma_1(s)\sigma(\bar{X}_{\eta(s)})).$$

We finally obtain the rate of convergence by applying once more the duality for stochastic integrals

$$\mathbb{E} f(X_T) - \mathbb{E} f(\bar{X}_T) = \mathbb{E} \int_0^T \int_{\eta(s)}^s D_u U_s^h \, du \, ds.$$

This last formula makes clear that $|\mathbb{E} f(X_T) - \mathbb{E} f(\bar{X}_T)| \leq T/n$ and leads to an expansion of $\mathbb{E} f(X_T) - \mathbb{E} f(\bar{X}_T)$ with some additional work (see Section 3).

In other stochastic equations, the error process also satisfies a linear equation [similar to (1)] but in a more general form; see Section 3.4. The aim of the next section is to establish a formula (called the duality formula) which will be a substitute for (2) when the error process $Y$ does not have an explicit expression.

### 3. Duality for the error process and application to delay equations.

3.1. *General form of the error process equation.* Throughout the paper, we consider a complete probability space $(\Omega, \mathcal{F}, \mathbb{P})$, which is the canonical space of a $d$-dimensional Brownian motion with finite horizon $W = \{(W_t^1, \ldots, W_t^d), 0 \leq t \leq T\}$. We denote by $\mathbf{F} = (\mathcal{F}_t)_{0 \leq t \leq T}$ the usual augmentation of the natural filtration of $W$. If $H$ is a separable Hilbert space and $p \in [1, +\infty)$, we denote by $L_a^p([0,T];H)$ the space of all measurable adapted processes $X = (X_t)_{0 \leq t \leq T}$, with values in $H$ such that $\mathbb{E} \int_0^T |X_t|^p \, dt < \infty$, where $|\cdot|$ denotes the norm on $H$.

Recall that according to the Itô representation theorem, any $H$-valued random variable $X$ such that $\mathbb{E}|X|^2 < \infty$ can be written in the following form:

$$X = \mathbb{E}(X) + \sum_{i=1}^d \int_0^T J_s^i(X) \, dW_s^i,$$

where $J^1(X), \ldots, J^d(X) \in L_a^2([0,T];H)$. Note that this defines $J^1, \ldots, J^d$ as linear operators mapping $L^2(\mathcal{F}_T; H)$ into $L_a^2([0,T];H)$. We will often use the notation $J(X)$ for the vector $(J^1(X), \ldots, J^d(X))$ and $Z \cdot dW$ for $\sum_{i=1}^d Z^i \, dW^i$.

We consider $\alpha_1, \ldots, \alpha_d, \beta, d+1$ linear continuous operators on the Hilbert space $L_a^2([0,T];H)$. The aim of Section 3.1 is to study the following generalization of (1):

$$(3) \qquad Y_t = \sum_{i=1}^d \int_0^t \alpha_i(Y)(s) \, dW_s^i + \int_0^t \beta(Y)(s) \, ds + G_t, \qquad 0 \leq t \leq T,$$

where $G \in L_a^2([0,T];H)$.

In fact, the study of (3) in the space $L_a^2$ will not be sufficient (see Section 3.4) and we will need $L^p$-estimates for $p$ large enough. For simplicity we state the assumption as follows.

ASSUMPTION A1.  For every $p \geq 2$, there exists a positive constant $C_p$ such that, for all $t \in [0, T]$, $Y \in L_a^p([0, T]; H)$,

$$(4) \quad \mathbb{E} \int_0^t \left( \sum_{i=1}^d |\alpha_i(Y)(s)|^2 \right)^{p/2} ds + \mathbb{E} \int_0^t |\beta(Y)(s)|^p \, ds \leq C_p \mathbb{E} \int_0^t |Y_s|^p \, ds.$$

PROPOSITION 1.  *Let Assumption* A1 *hold and let* $p \geq 2$. *Given* $G \in L_a^p([0, T]; H)$, (3) *has a unique solution* $Y \in L_a^p([0, T]; H)$ *and we have*

$$\mathbb{E} \int_0^T |Y_t|^p \, dt \leq \bar{C}_p \mathbb{E} \int_0^T |G_t|^p \, dt,$$

*for some constant* $\bar{C}_p$ *depending on* $C_p$ *(and* $T$) *only.*

PROOF.  Let $k$ be a positive real number. We define an equivalent norm $\| \cdot \|_k$ on the space $L_a^p([0, T]; H)$ by setting

$$\|Y\|_k^p = \mathbb{E} \int_0^T e^{-kt} |Y(t)|^p \, dt.$$

Let $A$ and $B$ be the operators defined by

$$(5) \qquad\qquad A(Y)(t) = \sum_{i=1}^d \int_0^t \alpha_i(Y)(s) \, dW_s^i,$$

$$(6) \qquad\qquad B(Y)(t) = \int_0^t \beta(Y)(s) \, ds.$$

We have, using the Burkholder–Davis–Gundy inequality (BDG inequality in the sequel) and $p \geq 2$,

$$\|A(Y)\|_k^p = \mathbb{E} \int_0^T dt \, e^{-kt} \left| \sum_{i=1}^d \int_0^t \alpha_i(Y)(s) \, dW_s^i \right|^p$$

$$\leq K_p \mathbb{E} \int_0^T dt \, e^{-kt} \left( \int_0^t \sum_{i=1}^d |\alpha_i(Y)(s)|^2 \, ds \right)^{p/2}$$

$$\leq K_p \mathbb{E} \int_0^T dt \, e^{-kt} t^{p/2-1} \int_0^t \left( \sum_{i=1}^d |\alpha_i(Y)(s)|^2 \right)^{p/2} ds$$

$$\leq K_p C_p T^{p/2-1} \mathbb{E} \int_0^T dt \, e^{-kt} \int_0^t ds \, |Y_s|^p,$$

where the constant $K_p$ comes from the BDG inequality and $C_p$ from Assumption A1. Using Fubini's theorem, we derive

$$(7) \qquad\qquad \|A(Y)\|_k^p \leq \frac{K_p C_p T^{p/2-1}}{k} \|Y\|_k^p.$$

A similar argument leads to

$$(8) \qquad \|B(Y)\|_k^p \leq \frac{C_p}{k} \|Y\|_k^p.$$

It follows that, for $k$ large enough, the operator $A+B$ is a contraction for the norm $\|\cdot\|_k$. Hence, if $I$ denotes the identity, the operator $I-(A+B)$ (acting on $L_a^p([0,T];H)$) is invertible and this implies existence and uniqueness for the solution of (3). $\square$

REMARK 2. Note that if the process $G$ has right-continuous (resp. continuous) paths, the solution of (3) has a right-continuous (resp. continuous) modification. This modification will still be denoted by $Y$.

3.2. *The duality formulae.* The purpose of this section is to establish two duality formulae relating $\mathbb{E}\langle F, Y_T\rangle$ to $G_T$ for $F$ in $L^2(\Omega)$. We first introduce some notation. We denote by $\langle\cdot,\cdot\rangle$ the inner product on $H$. The operators $A^*: L_a^2([0,T];H) \to L_a^2([0,T];H)$ and $B^*: L_a^2([0,T];H) \to L_a^2([0,T];H)$ are the adjoints of the operators $A$ and $B$ defined by (5) and (6). If we set $\alpha(Y)=(\alpha_1(Y),\dots,\alpha_d(Y))$, we can view $\alpha$ as a linear operator mapping the Hilbert space $L_a^2([0,T];H)$ into $L_a^2([0,T];H^d)$. The adjoint operator $\alpha^*$ maps $L_a^2([0,T];H^d)$ into $L_a^2([0,T];H)$. The following proposition relates the operators $A^*$ and $B^*$ to $\alpha^*$ and $\beta^*$.

PROPOSITION 3. *The operators $A^*$ and $B^*$ are given by*

$$A^*(Z)(t) = \sum_{i=1}^d \alpha_i^* \left( J^i \left( \int_0^T Z_s\,ds \right) \right)(t)$$

$$= \alpha^* \left( J \left( \int_0^T Z_s\,ds \right) \right)(t),$$

$$B^*(Z)(t) = \beta^* \left( \mathbb{E} \left( \int_\cdot^T Z_s\,ds \Big| \mathcal{F}_\cdot \right) \right)(t).$$

PROOF. For all $Y, Z \in L_a^2([0,T];H)$, we have

$$\mathbb{E}\int_0^T \langle Z_t, A(Y)(t)\rangle\,dt = \mathbb{E}\sum_{i=1}^d \int_0^T \left\langle Z_t, \int_0^t \alpha_i(Y)(s)\,dW_s^i \right\rangle dt$$

$$= \sum_{i=1}^d \int_0^T \mathbb{E}\int_0^t \langle J_s^i(Z_t), \alpha_i(Y)(s)\rangle\,ds\,dt$$

$$= \mathbb{E}\int_0^T \left( \int_s^T \sum_{i=1}^d \langle J_s^i(Z_t), \alpha_i(Y)(s)\rangle\,dt \right) ds.$$

Note that, since $Z_t$ is $\mathcal{F}_t$-measurable, we have $J_s^i(Z_t) = 0$ for $t < s$. Hence, using the linearity of the operators $J^i$,

$$\mathbb{E}\int_0^T \langle Z_t, A(Y)(t)\rangle\, dt = \mathbb{E}\int_0^T\left(\int_0^T\sum_{i=1}^d \langle J_s^i(Z_t), \alpha_i(Y)(s)\rangle\, dt\right) ds$$

$$= \mathbb{E}\int_0^T \sum_{i=1}^d \left\langle J_s^i\left(\int_0^T Z_t\, dt\right), \alpha_i(Y)(s)\right\rangle ds$$

$$= \mathbb{E}\int_0^T \left\langle \sum_{i=1}^d \alpha_i^*\left(J^i\left(\int_0^T Z_t\, dt\right)\right)(s), Y(s)\right\rangle ds,$$

which proves the formula for $A^*$. We proceed similarly with $B$:

$$\mathbb{E}\int_0^T \langle Z_t, B(Y)(t)\rangle\, dt = \mathbb{E}\int_0^T\left\langle Z_t, \int_0^t \beta(Y)(s)\, ds\right\rangle dt$$

$$= \mathbb{E}\int_0^T\left\langle \int_s^T Z_t\, dt, \beta(Y)(s)\right\rangle ds$$

$$= \mathbb{E}\int_0^T\left\langle \beta^*\left(\mathbb{E}\left(\int_{\cdot}^T Z_t\, dt\Big|\mathcal{F}.\right)\right)(s), Y(s)\right\rangle ds. \quad \square$$

The next theorem states the two basic duality formulae. In order to relate $\mathbb{E}\langle\Phi, Y_T\rangle$ to $G_T$, we need a formula for $\mathbb{E}\int_0^T\langle F_t, Y_t\rangle\, dt$ when $(F_t)_{0\le t\le T} \in L_a^2([0,T]; H)$.

THEOREM 4.    Let $G \in L_a^2([0,T]; H)$ and let $Y$ be the solution of (3).

1. If $F = (F_t)_{0\le t\le T} \in L_a^2([0,T]; H)$, we have

$$\mathbb{E}\int_0^T\langle F_t, Y_t\rangle\, dt = \mathbb{E}\int_0^T\langle\theta_t, G_t\rangle\, dt,$$

with $\theta = (I - A^* - B^*)^{-1}(F)$.

2. If $G$ has a continuous modification, then $Y$ has a continuous modification (which we still denote by $Y$), and if $\Phi$ is an $\mathcal{F}_T$-measurable square integrable random variable with values in $H$, we have

$$\mathbb{E}\langle\Phi, Y_T\rangle = \mathbb{E}\langle\Phi, G_T\rangle + \mathbb{E}\int_0^T\langle\hat{\theta}_t, G_t\rangle\, dt,$$

with

$$\hat{\theta} = (I - A^* - B^*)^{-1}[\alpha^*(J(\Phi)) + \beta^*(\mathbb{E}(\Phi|\mathcal{F}.))].$$

PROOF.  The first part of the theorem comes from the equality $Y = (I - A - B)^{-1}(G)$ and standard duality theory. For the second part, it is clear from (3) that if $t \mapsto G_t$ is continuous, the process $Y$ has a continuous modification. We also have

$$\mathbb{E}\langle \Phi, Y_T \rangle = \mathbb{E}\langle \Phi, G_T \rangle + \mathbb{E}\left\langle \Phi, \int_0^T \alpha(Y)(s) \cdot dW_s + \int_0^T \beta(Y)(s)\, ds \right\rangle.$$

Now

$$\mathbb{E}\left\langle \Phi, \int_0^T \alpha(Y)(s) \cdot dW_s \right\rangle = \mathbb{E}\int_0^T \sum_{i=1}^d \langle J_s^i(\Phi), \alpha_i(Y)(s) \rangle\, ds$$

and

$$\mathbb{E}\left\langle \Phi, \int_0^T \beta(Y)(s)\, ds \right\rangle = \mathbb{E}\int_0^T \langle \mathbb{E}(\Phi|\mathcal{F}_s), \beta(Y)(s) \rangle\, ds.$$

Therefore

$$\mathbb{E}\langle \Phi, Y_T \rangle = \mathbb{E}\langle \Phi, G_T \rangle + \mathbb{E}\int_0^T \langle F_s, Y_s \rangle\, ds,$$

with $F_t = \alpha^*(J(\Phi))(t) + \beta^*(\mathbb{E}(\Phi|\mathcal{F}_.))(t)$. And the result follows from the first part. □

REMARK 5. The processes $\theta$ and $\hat{\theta}$ given by Theorem 4 can also be characterized in connection with backward stochastic differential equations (BSDEs). First note that $\theta$ satisfies the dual equation

$$(9) \qquad \theta = F + A^*(\theta) + B^*(\theta),$$

which, using Proposition 3, can be written

$$(10) \qquad \theta = F + \alpha^*\left(J\left(\int_0^T \theta_s\, ds\right)\right) + \beta^*\left(\mathbb{E}\left(\int_.^T \theta_s\, ds \Big| \mathcal{F}_.\right)\right).$$

Now, observe that, given $\theta \in L_a^2([0,T];H)$ and $\Psi \in L^2(\mathcal{F}_T;H)$, the pair of processes $(\tilde{Y}, \tilde{Z})$, defined by

$$\tilde{Y}_t = \mathbb{E}\left(\Psi + \int_t^T \theta_s\, ds \Big| \mathcal{F}_t\right) \quad \text{and} \quad \tilde{Z}_t = J_t\left(\Psi + \int_0^T \theta_s\, ds\right),$$

is the unique solution of the following BSDE:

$$d\tilde{Y}_t = -\theta_t\, dt + \tilde{Z}_t \cdot dW_t,$$
$$\tilde{Y}_T = \Psi.$$

It follows that, if $\theta$ satisfies (10), we have

$$\theta_t = F_t + \alpha^*(Z)(t) + \beta^*(Y)(t),$$

where the pair $(Y, Z)$ solves the following BSDE:

$$dY_t = -(F_t + \alpha^*(Z)(t) + \beta^*(Y)(t)) \, dt + Z_t \cdot dW_t,$$
$$Y_T = 0.$$

Similarly, the process $\hat{\theta}$ in Theorem 4 is given by

$$\hat{\theta}_t = \alpha^*(\hat{Z})(t) + \beta^*(\hat{Y})(t),$$

where the pair $(\hat{Y}, \hat{Z})$ solves the following BSDE:

$$d\hat{Y}_t = -(\alpha^*(\hat{Z})(t) + \beta^*(\hat{Y})(t)) \, dt + \hat{Z}_t \cdot dW_t,$$
$$\hat{Y}_T = \Phi.$$

3.3. *Estimates for the stochastic derivatives of the dual equation.* In this section we study the dual equation (9). We will establish estimates for the derivatives of the process $\theta$ involved in the duality formulae. These estimates will be useful for the study of the weak rate of convergence of the Euler scheme for the computation of expectations of regular functions (see Section 3.4). Under Assumption A1, we know from Proposition 1 that the operator $I - A - B$ is invertible when considered on the space $L_a^p([0, T]; H)$ for $p \geq 2$. The operators $A^*$ and $B^*$ can be viewed as bounded linear operators on $L_a^p([0, T]; H)$ for $p \leq 2$, and it follows that the operator $(I - A^* - B^*)$ is invertible on the space $L_a^p([0, T]; H)$, for $p \leq 2$.

In other words, we may assert that, given $F \in L_a^p([0, T]; H)$ (with $1 < p \leq 2$), there exists a unique $\theta \in L_a^p([0, T]; H)$, satisfying

$$\theta = F + A^*(\theta) + B^*(\theta).$$

Using Proposition 3, we have

$$\tag{11} \theta = F + \alpha^* \left( J \left( \int_0^T \theta_s \, ds \right) \right) + \beta^* \left( \mathbb{E} \left( \int_{\cdot}^T \theta_s \Big| \mathcal{F}_{\cdot} \right) \right).$$

We want to differentiate this equation, in order to estimate the Malliavin derivatives of $\theta$. We will need some regularity assumptions on the operators $\alpha$ and $\beta$. For the basic theory of Sobolev spaces on Wiener space and for standard notation, we refer the reader to [13].

We denote by $D$ the derivative operator. If $X$ is a simple functional with values in $H$, $DX$ is a random variable with values in $L^2([0, T]; H^d)$ and can be viewed as a nonadapted process $(D_t X)_{0 \leq t \leq T}$. We will say that a functional (or random variable) $X$ with values in $H$ is smooth, if it can be written as a finite sum of multiple Wiener integrals with continuous deterministic integrands. Note that if $X$ is smooth, the process $(D_t X)$ has a right-continuous modification. We denote by $\mathcal{S}_a([0, T]; H)$ the space of all adapted processes $(Y_t)_{0 \leq t \leq T}$ with values in $H$, with continuous paths such that, for each $t \in [0, T]$, the random variable $Y_t$ is smooth in the previous sense. The space $\mathcal{S}_a([0, T]; H)$ is dense in $L_a^p([0, T]; H)$ for all $p \in [1, +\infty)$. Our regularity assumptions on $\alpha$ and $\beta$ can now be formulated as follows.

ASSUMPTION A2. For $\gamma = \alpha$ and $\gamma = \beta$ and for all $Z \in \mathcal{S}_a([0, T]; H)$, the process $\gamma^*(Z)$ is in $\mathbb{D}^{2,2}$ and, for all $u, v \in [0, T]$, there exist operators denoted by $D_u\gamma^*$, $D_{uv}^2\gamma^*$, such that, for $Z \in \mathcal{S}_a([0, T]; H)$,

$$D_u(\gamma^*(Z)) = (D_u\gamma^*)(Z) + \gamma^*(D_u Z),$$
$$D_{uv}^2(\gamma^*(Z)) = (D_{uv}^2\gamma^*)(Z) + (D_u\gamma^*)(D_v Z)$$
$$+ (D_v\gamma^*)(D_u Z) + \gamma^*(D_{uv}^2 Z).$$

Moreover, these operators satisfy the following estimate, for all $q_1, q_2$ with $1 \leq q_1 < q_2 \leq 2$:

$$
\begin{aligned}
(12) \quad &\left( \mathbb{E} \int_0^T |D_u\gamma^*(Z)(t)|^{q_1} + |D_{uv}^2\gamma^*(Z)(t)|^{q_1} \, dt \right)^{1/q_1} \\
&\leq C_{q_1,q_2} \left( \mathbb{E} \int_0^T |Z_t|^{q_2} \, dt \right)^{1/q_2},
\end{aligned}
$$

where the constants $C_{q_1,q_2}$ do not depend on $u, v$ (or $Z$).

The estimate (12) for $\alpha$ and $\beta$ allows us to extend the operators $D_u\alpha^*$, $D_{uv}^2\alpha^*$, $D_u\beta^*$, $D_{uv}^2\beta^*$ as continuous operators from $L_a^{q_2}([0, T]; H)$ into $L_a^{q_1}([0, T]; H)$, for $1 \leq q_1 < q_2 \leq 2$. Note that in typical examples [see Section 3.4, (21), (22)], the operators $D_u\alpha^*$, $D_u\beta^*$ are not bounded from $L_a^p([0, T]; H)$ into itself. In the sequel, it will be convenient to use the following notation. For a random variable $Z$ in $\mathbb{D}^{2,2}$, and $u, v \in [0, T]$, let

$$n(Z, u, v) = |Z| + |D_u Z| + |D_v Z| + |D_{uv}^2 Z|.$$

PROPOSITION 6. Let Assumptions A1 and A2 hold. Let $F \in L_a^2([0, T]; H)$. We assume that $F$ has a modification $F_t$ such that, for each $t \in [0, T]$, $F_t \in \mathbb{D}^{2,2}$ and $\sup_{0 \leq u, v \leq T} \mathbb{E} \int_0^T n(F_t, u, v)^2 \, dt < \infty$. Then, the solution of (11) has a modification $\theta_t$ satisfying, for $1 \leq p < q \leq 2$,

$$\left( \mathbb{E} \int_0^T n^p(\theta_t, u, v) \, dt \right)^{1/p} \leq K_{p,q} \left( \mathbb{E} \int_0^T n^q(F_t, u, v) \, dt \right)^{1/q},$$

where the constants $K_{p,q}$ depend on $T$ and the constants $C_p$ and $C_{q_1,q_2}$ in Assumptions A1 and A2 (and not on $u, v$ or $F$).

For the proof of Proposition 6, we will need the following commutation relations between the operators $J$ and the derivative and conditional expectation operators. We omit the proof which involves only classical arguments of analysis on Wiener space.

LEMMA 7.   *For all $X \in \mathbb{D}^{1,2}(H)$, we have*

$$D_u(J(X)(v)) = J(D_u(X))(v)\mathbb{1}_{\{u \leq v\}} \quad and$$

$$D_u(\mathbb{E}(X|\mathcal{F}_v)) = \mathbb{E}(D_u(X)|\mathcal{F}_v)\mathbb{1}_{\{u \leq v\}}, \qquad du\, dv \ \ a.e.$$

PROOF OF PROPOSITION 6.   The formal differentiation of (11) gives

$$D_u\theta_t = D_uF_t + (D_u\alpha^*)\left(J\left(\int_0^T \theta_s\, ds\right)\right)(t) + (D_u\beta^*)\left(\mathbb{E}\left(\int_{\cdot}^T \theta_s\, ds\Big|\mathcal{F}_{\cdot}\right)\right)(t)$$

$$+ \alpha^*\left(D_u\left(J\left(\int_0^T \theta_s\, ds\right)\right)\right)(t) + \beta^*\left(D_u\left(\mathbb{E}\left(\int_{\cdot}^T \theta_s\, ds\Big|\mathcal{F}_{\cdot}\right)\right)\right)(t).$$

Using Lemma 7 and the linearity of $D_u$, we get

$$
\begin{aligned}
D_u\theta_t = {}& D_uF_t + (D_u\alpha^*)\left(J\left(\int_0^T \theta_s\, ds\right)\right)(t) \\
& + (D_u\beta^*)\left(\mathbb{E}\left(\int_{\cdot}^T \theta_s\, ds\Big|F_{\cdot}\right)\right)(t) \\
& + \alpha^*\left(\mathbb{1}_{[u,T]}J\left(\int_0^T D_u\theta_s\, ds\right)\right)(t) \\
& + \beta^*\left(\mathbb{1}_{[u,T]}\mathbb{E}\left(\int_{\cdot}^T D_u\theta_s\, ds\Big|\mathcal{F}_{\cdot}\right)\right)(t).
\end{aligned}
$$

(13)

Now let $I_{u,T}$ be the operator on $L_a^2([0,T];H)$ defined by

$$I_{u,T}(Y) = \mathbb{1}_{[u,T]}Y.$$

Clearly, $I_{u,T}$ is a self-adjoint operator and defines an operator with norm 1 on $L_a^p([0,T];H)$ for every $p \in [1,+\infty)$. We have

$$\alpha^* \circ I_{u,T} = (I_{u,T} \circ \alpha)^* \quad \text{and} \quad \beta^* \circ I_{u,T} = (I_{u,T} \circ \beta)^*,$$

and the operators $I_{u,T} \circ \alpha$ and $I_{u,T} \circ \beta$ satisfy Assumption A1 with the same constants $C_p$ as $\alpha$ and $\beta$. Now let $A_u$ and $B_u$ be the operators defined on $L_a^p([0,T];H)$ $(p \geq 2)$ by

$$A_u(Y)(t) = \sum_{i=1}^d \int_0^t \mathbb{1}_{[u,T]}(s)\alpha_i(Y)(s)\, dW_s^i,$$

$$B_u(Y)(t) = \int_0^t \mathbb{1}_{[u,T]}(s)\beta(Y)(s)\, ds.$$

Using Proposition 3, we can rewrite (13) as

$$D_u\theta = \phi_u + (A_u^* + B_u^*)(D_u\theta),$$

where $\phi_u$ is the adapted process defined by

$$\phi_u(t) = D_u F_t + (D_u \alpha^*)\left(J\left(\int_0^T \theta_s\, ds\right)\right)(t)$$

$$+ (D_u \beta^*)\left(\mathbb{E}\left(\int_\cdot^T \theta_s\, ds \Big| \mathcal{F}_\cdot\right)\right)(t).$$

It follows from Proposition 1 and the properties of the operators $I_{u,T} \circ \alpha$ and $I_{u,T} \circ \beta$ that, for $p \in [1, 2]$,

$$\mathbb{E}\int_0^T |D_u \theta_t|^p\, dt \le C_p \mathbb{E}\int_0^T |\phi_u(t)|^p\, dt,$$

where $C_p$ does not depend on $u$. On the other hand, we deduce from Assumption A2 and the assumptions on $F$ that, for $1 \le p < q \le 2$,

$$\left(\mathbb{E}\int_0^T |\phi_u(t)|^p\, dt\right)^{1/p} \le C_p\left(\mathbb{E}\int_0^T |D_u F_t|^p\, dt\right)^{1/p}$$

$$+ C_{p,q}\left\{\left(\mathbb{E}\int_0^T \Big|J_t\left(\int_0^T \theta_s\, ds\right)\Big|^q\, dt\right)^{1/q}\right.$$

$$\left. + \left(\mathbb{E}\int_0^T \Big|\mathbb{E}\left(\int_t^T \theta_s\, ds \Big| \mathcal{F}_t\right)\Big|^q\, dt\right)^{1/q}\right\},$$

where here again the constants $C_p$ and $C_{p,q}$ do not depend on $u$. Here we have used the notation $J_t(Z)$ for $J(Z)(t)$. We have, using Jensen's and Hölder's inequalities,

$$\left(\mathbb{E}\int_0^T \Big|\mathbb{E}\left(\int_t^T \theta_s\, ds \Big| \mathcal{F}_t\right)\Big|^q\, dt\right)^{1/q} \le \left(\mathbb{E}\int_0^T \Big|\int_t^T \theta_s\, ds\Big|^q\, dt\right)^{1/q}$$

$$\le T\left(\mathbb{E}\int_0^T |\theta_t|^q\, dt\right)^{1/q}.$$

On the other hand, using $1 < q \le 2$ and the BDG inequality, we have

$$\left(\mathbb{E}\int_0^T \Big|J_t\left(\int_0^T \theta_s\, ds\right)\Big|^q\, dt\right)^{1/q}$$

$$\le T^{1/q - 1/2}\left(\mathbb{E}\left(\int_0^T \Big|J_t\left(\int_0^T \theta_s\, ds\right)\Big|^2\, dt\right)^{q/2}\right)^{1/q}$$

$$\le T^{1/q - 1/2} C_q\left(\mathbb{E}\Big|\int_0^T \theta_s\, ds\Big|^q\right)^{1/q}$$

$$\le T^{1/2} C_q\left(\mathbb{E}\int_0^T |\theta_s|^q\, ds\right)^{1/q}.$$

Recall that $\theta = (I - A^* - B^*)^{-1}(F)$, so that

$$\left( \mathbb{E} \int_0^T |\theta_t|^q \, dt \right)^{1/q} \leq C_q \left( \mathbb{E} \int_0^T |F_t|^q \, dt \right)^{1/q}.$$

Hence, we have, for $1 \leq p < q \leq 2$,

$$\left( \mathbb{E} \int_0^T |D_u \theta_t|^p \, dt \right)^{1/p} \leq C_p \left( \mathbb{E} \int_0^T |D_u F_t|^p \, dt \right)^{1/p} + C_{p,q} \left( \mathbb{E} \int_0^T |F_t|^q \, dt \right)^{1/q},$$

where the constants $C_p$ and $C_{p,q}$ depend on $T$ but not on $u$. We can now differentiate (13) with respect to $v$ and derive in a similar manner the estimate

$$\left( \mathbb{E} \int_0^T |D_{uv}^2 \theta_t|^p \, dt \right)^{1/p} \leq C_p \left( \mathbb{E} \int_0^T |D_{uv}^2 F_t|^p \, dt \right)^{1/p}$$

$$+ C_{p,q} \left( \mathbb{E} \int_0^T (|F_t| + |D_u F_t| + |D_v F_t|)^q \, dt \right)^{1/q}.$$

Note that in order to justify the formal differentiations, it suffices to use an iterating procedure of the form $\theta_0 = F$ and $\theta_{j+1} = F + (A^* + B^*)(\theta_j)$. The regularity of $\theta_j$ carries over to $\theta_{j+1}$ and we have convergence of $\theta_j$ toward $\theta$, together with the derivatives, with the suitable norms. $\square$

For the application of our methodology to numerical schemes, we will need to consider families of operators $(\alpha^h, \beta^h)_{h \geq 0}$ and introduce some additional notation. Let $\mathcal{A}_{1,2}$ denote the space of all operators $(\alpha, \beta)$ satisfying Assumptions A1 and A2. For $(\alpha, \beta) \in \mathcal{A}_{1,2}$, let $C_p(\alpha, \beta)$ [resp. $C_{q_1,q_2}(\alpha, \beta)$] be the smallest constant for which (4) [resp. (12)] holds. The following proposition follows easily from Proposition 6.

PROPOSITION 8. *Assume $(\alpha^h, \beta^h)_{h \geq 0}$ is a family of operators in $\mathcal{A}_{1,2}$, satisfying, for all $p \in [2, +\infty)$ and for all $q_1, q_2$ with $1 \leq q_1 < q_2 \leq 2$,*

$$\lim_{h \to 0} C_p(\alpha^h - \alpha, \beta^h - \beta) = 0 \quad \text{and} \quad \lim_{h \to 0} C_{q_1, q_2}(\alpha^h - \alpha, \beta^h - \beta) = 0.$$

*If $(F^h)_{h \geq 0}$ is a family of adapted processes satisfying $\sup_{0 \leq u, v \leq T} \mathbb{E} \int_0^T n(F_t^h, u, v)^2 \, dt < \infty$ and*

$$\lim_{h \to 0} \sup_{0 \leq u, v \leq T} \mathbb{E} \int_0^T n(F_t^h - F_t^0, u, v)^2 \, dt = 0,$$

*then the processes $\theta^h$ [defined by $\theta^h = (I - A_h^* - B_h^*)^{-1}(F^h)$] satisfy*

$$\forall p \in [1, 2) \qquad \lim_{h \to 0} \sup_{0 \leq u, v \leq T} \mathbb{E} \int_0^T n(\theta_t^h - \theta_t^0, u, v)^p \, dt = 0.$$

3.4. *Application to the Euler approximation of delay equations*: *the regular case.* In this section we assume to simplify the notation that all processes take values in $\mathbb{R}$. We are interested in the expansion of $\mathbb{E}f(X_T) - \mathbb{E}f(\bar{X}_T)$ where the process $(X_t)_{t \in [-r,T]}$ solves the stochastic delay equation

$$dX_t = \sigma\left(\int_{-r}^0 X_{t+s}\, d\nu(s)\right) dW_t + b\left(\int_{-r}^0 X_{t+s}\, d\nu(s)\right) dt, \qquad t \geq 0,$$

$$X_s = \xi_s, \qquad s \in [-r, 0],$$

where $r > 0$, $\xi \in \mathcal{C}^1([-r,0], \mathbb{R})$ and $\nu$ is a finite measure.

We denote by $\bar{X}$ the following Euler approximation of $(X_t)_{t \in [-r,T]}$ with step $h = r/n$:

$$d\bar{X}_t = \sigma\left(\int_{-r}^0 \bar{X}_{\eta(t)+\eta(s)}\, d\nu(s)\right) dW_t + b\left(\int_{-r}^0 \bar{X}_{\eta(t)+\eta(s)}\, d\nu(s)\right) dt, \qquad t \geq 0,$$

$$\bar{X}_s = \xi_s, \qquad s \in [-r, 0],$$

with $\eta(s) = \frac{[ns/r]}{n/r}$, where $[t]$ stands for the entire part of $t$. We assume that the functions, $f$, $\sigma$ and $b$ are in $\mathcal{C}_b^3$. Note that if $\nu$ is the Dirac measure at $0$, we have a standard diffusion.

For existence, uniqueness and moment estimates of solutions of stochastic delay equations in the above form, we refer to [12]. Consistency of the Euler scheme for delay equations in the case where $\nu$ is a Dirac mass has been studied in [9] through an extension of the PDE method. An infinite-dimensional extension of the PDE method is used in [4] but their result is limited to drift coefficients linearly dependent on the past and nondelayed diffusion coefficient.

Our expansion result is derived from the duality formula and the following lemma.

LEMMA 9. *Let $(U_s^h)$ be a family of $\mathcal{F}_T$-measurable real-valued processes, such that $\forall s \in [0, T]$, $\forall h \in [0,1]$, $U_s^h \in \mathbb{D}^{1,2}$. We assume that, for some $p > 1$, we have*

$$\lim_{h \to 0} \mathbb{E}\int_0^T |U_s^h - U_s^0|^p\, ds = 0$$

*and*

$$\lim_{h \to 0}\int_{-r}^0\int_0^T \sup_{v \in [\eta(s)+\eta(u), s+u]} \|D_v U_s^h - D_{s+u}U_s^0\|_p \mathbb{1}_{\eta(s)+\eta(u) \geq 0}\, ds\, d\nu(u) = 0.$$

*Then*

$$\mathbb{E}\int_0^T U_s^h \int_{-r}^0 (\bar{X}_{s+u} - \bar{X}_{\eta(s)+\eta(u)})\, d\nu(u)\, ds = hC(U^0) + I^h(U^0) + o(h),$$

*where*

$$C(U^0) = \tfrac{1}{2} \int_{-r}^0 \mathbb{E} \int_0^T [U_s^0(\tilde{b}_{s+u}^0 \mathbb{1}_{s+u \geq 0} + \xi'_{s+u} \mathbb{1}_{s+u < 0})$$
$$+ \tilde{\sigma}_{s+u}^0 D_{s+u} U_s^0 \mathbb{1}_{s+u \geq 0}] \, ds \, d\nu(u),$$

$$I^h(U^0) = \int_{-r}^0 \mathbb{E} \int_0^T [U_s^0(\tilde{b}_{s+u}^0 \mathbb{1}_{s+u \geq 0} + \xi'_{s+u} \mathbb{1}_{s+u < 0})$$
$$+ \tilde{\sigma}_{s+u}^0 D_{s+u} U_s^0 \mathbb{1}_{s+u \geq 0}] \, ds(u - \eta(u)) \, d\nu(u),$$

$$\tilde{b}_t^0 = b\left(\int_{-r}^0 X_{t+v} \, d\nu(v)\right) \quad and \quad \tilde{\sigma}_t^0 = \sigma\left(\int_{-r}^0 X_{t+v} \, d\nu(v)\right).$$

REMARK. Observe that $|I^h(U^0)| \leq hC$. Moreover, when $\nu$ is a Dirac mass at $-r$ a good discretization of the time interval gives $I^h(U^0) = 0$. When $\nu$ is an absolutely continuous measure, we obtain $I^h(U^0) = hC(U^0) + o(h)$.

PROOF OF LEMMA 9. We have

$$\bar{X}_{s+u} - \bar{X}_{\eta(s)+\eta(u)} = \left(\int_{\eta(s)+\eta(u)}^{s+u} \tilde{\sigma}_{\eta(t)}^h \, dW_t + \int_{\eta(s)+\eta(u)}^{s+u} \tilde{b}_{\eta(t)}^h \, dt\right) \mathbb{1}_{\eta(s)+\eta(u) \geq 0}$$
$$+ (\xi_{s+u} - \xi_{\eta(s)+\eta(u)}) \mathbb{1}_{s+u < 0}$$
$$+ (\xi_0 - \xi_{\eta(s)+\eta(u)}) \mathbb{1}_{\eta(s)+\eta(u) \leq 0 < s+u},$$

where

$$(14) \qquad\qquad \tilde{\sigma}_t^h = \sigma\left(\int_{-r}^0 \bar{X}_{t+\eta(v)} \, d\nu(v)\right)$$

and

$$(15) \qquad\qquad \tilde{b}_t^h = b\left(\int_{-r}^0 \bar{X}_{t+\eta(v)} \, d\nu(v)\right).$$

This gives

$$\mathbb{E} \int_0^T U_s^h \int_{-r}^0 (\bar{X}_{s+u} - \bar{X}_{\eta(s)+\eta(u)}) \, d\nu(u) \, ds$$
$$= \int_{-r}^0 \mathbb{E} \int_0^T U_s^h \mathbb{1}_{\eta(s)+\eta(u) \geq 0} \int_{\eta(s)+\eta(u)}^{s+u} \tilde{\sigma}_{\eta(t)}^h \, dW_t \, ds \, d\nu(u)$$
$$+ \int_{-r}^0 \mathbb{E} \int_0^T U_s^h \mathbb{1}_{\eta(s)+\eta(u) \geq 0} \int_{\eta(s)+\eta(u)}^{s+u} \tilde{b}_{\eta(t)}^h \, dt \, ds \, d\nu(u)$$
$$+ \int_{-r}^0 \mathbb{E} \int_0^T U_s^h (\xi_{s+u} - \xi_{\eta(s)+\eta(u)}) \mathbb{1}_{s+u < 0} \, ds \, d\nu(u)$$
$$+ \int_{-r}^0 \mathbb{E} \int_0^T U_s^h (\xi_0 - \xi_{\eta(s)+\eta(u)}) \mathbb{1}_{\eta(s)+\eta(u) \leq 0 < s+u} \, ds \, d\nu(u).$$

By duality, we obtain for the first term

$$\int_{-r}^0 \mathbb{E} \int_0^T U_s^h \mathbb{1}_{\eta(s)+\eta(u)\geq 0} \int_{\eta(s)+\eta(u)}^{s+u} \tilde{\sigma}_{\eta(t)}^h \, dW_t \, ds \, d\nu(u)$$

$$= \int_{-r}^0 \mathbb{E} \int_0^T \mathbb{1}_{\eta(s)+\eta(u)\geq 0} \int_{\eta(s)+\eta(u)}^{s+u} D_t U_s^h \tilde{\sigma}_{\eta(t)}^h \, dt \, ds \, d\nu(u).$$

Since $\eta(s)+\eta(u) \leq 0 < s+u$ implies $-u < s < -u+2h$, one can easily verify that

$$\int_{-r}^0 \mathbb{E} \int_0^T U_s^h (\xi_0 - \xi_{\eta(s)+\eta(u)}) \mathbb{1}_{\eta(s)+\eta(u)\leq 0 < s+u} \, ds \, d\nu(u) = o(h).$$

Now recall that if $g$ is an integrable function on $[0,T]$, we have $\int_0^T g(s)(s - \eta(s)) \, ds = \frac{h}{2} \int_0^T g(s) \, ds + o(h)$. It follows that in order to derive the expansion stated in the lemma, it is enough to prove that when $h$ tends to 0,

$$(16) \quad \begin{aligned} &\frac{1}{h} \Big| \int_{-r}^0 \mathbb{E} \int_0^T \overline{U}_s^h \int_{\eta(s)+\eta(u)}^{s+u} \tilde{b}_{\eta(t)}^h \, dt \, ds \, d\nu(u) \\ &\qquad - \int_{-r}^0 \mathbb{E} \int_0^T \overline{U}_s^0 \tilde{b}_{s+u}^0 \Delta_{s+u} \, ds \, d\nu(u) \Big| \to 0, \end{aligned}$$

$$(17) \quad \begin{aligned} &\frac{1}{h} \Big| \int_{-r}^0 \mathbb{E} \int_0^T \int_{\eta(s)+\eta(u)}^{s+u} D_t \overline{U}_s^h \tilde{\sigma}_{\eta(t)}^h \, dt \, ds \, d\nu(u) \\ &\qquad - \int_{-r}^0 \mathbb{E} \int_0^T \tilde{\sigma}_{s+u}^0 D_{s+u} \overline{U}_s^0 \Delta_{s+u} \, ds \, d\nu(u) \Big| \to 0, \end{aligned}$$

$$(18) \quad \begin{aligned} &\frac{1}{h} \Big| \int_{-r}^0 \mathbb{E} \int_0^T U_s^h \int_{\eta(s)+\eta(u)}^{s+u} \xi_t' \, dt \, \mathbb{1}_{s+u<0} \, ds \, d\nu(u) \\ &\qquad - \int_{-r}^0 \mathbb{E} \int_0^T U_s^0 \xi_{s+u}' \mathbb{1}_{s+u<0} \Delta_{s+u} \, ds \, d\nu(u) \Big| \to 0, \end{aligned}$$

where $\Delta_{s+u} = s+u-\eta(s)-\eta(u)$, $\overline{U}_s^h = U_s^h \mathbb{1}_{\eta(s)+\eta(u)\geq 0}$ and $\overline{U}_s^0 = U_s^0 \mathbb{1}_{s+u\geq 0}$. Using Hölder's inequality, the left-hand side of (16) is bounded by

$$\int_0^T \|U_s^h - U_s^0\|_p \, ds \sup_t \|\tilde{b}_t^h\|_q$$

$$+ \int_0^T \|U_s^0\|_p \, ds \sup_{u,s} \sup_{t\in[\eta(s)+\eta(u),s+u]} \|\tilde{b}_{\eta(t)}^h - \tilde{b}_{s+u}^0\|_q$$

$$+ \sup_u \int_0^T \|U_s^h\|_p \mathbb{1}_{\eta(s)+\eta(u)\leq 0\leq s+u} \, ds \sup_t \|\tilde{b}_t^h\|_q,$$

with $\frac{1}{p} + \frac{1}{q} = 1$. Observe that $\sup_u \int_0^T \|U_s^h\|_p \mathbb{1}_{\eta(s)+\eta(u)\le 0 \le s+u}\, ds \sup_t \|\tilde{b}_t^h\|_q = o(h)$. In the same way the left-hand side of (17) is bounded by

$$\int_{-r}^0 \int_0^T \sup_{t\in[\eta(s)+\eta(u),s+u]} \|D_t U_s^h - D_{s+u} U_s^0\|_p \mathbb{1}_{\eta(s)+\eta(u)\ge 0}\, ds\, d\nu(u) \sup_t \|\tilde{\sigma}_t^h\|_q$$

$$+ \int_{-r}^0 \int_0^T \|D_{s+u} U_s^0\|_p\, ds\, d\nu(u) \sup_{u,s} \sup_{t\in[\eta(s)+\eta(u),s+u]} \|\tilde{\sigma}_{\eta(t)}^h - \tilde{\sigma}_{s+u}^0\|_q$$

$$+ o(h).$$

Similar bounds hold for the left-hand side of (18) and we obtain the result of Lemma 9 as soon as

$$\sup_{h>0, s\in[0,T]} \|\tilde{b}_s^h\|_q < \infty \quad \text{and} \quad \sup_{h>0, s\in[0,T]} \|\tilde{\sigma}_s^h\|_q < \infty,$$

$$\sup_{u,s} \sup_{t\in[\eta(s)+\eta(u),s+u]} \|\tilde{\sigma}_{\eta(t)}^h - \tilde{\sigma}_{s+u}^0\|_q \to 0$$

and

$$\sup_{u,s} \sup_{t\in[\eta(s)+\eta(u),s+u]} \|\tilde{b}_{\eta(t)}^h - \tilde{b}_{s+u}^0\|_q \to 0,$$

for all $q \in [1, +\infty)$. This can be proved as in [6].   $\square$

The following theorem can be viewed as an analogue of the classical expansion of the error for diffusions (see [15]).

THEOREM 10.   *We have, if $f \in \mathcal{C}_b^3$,*

$$\mathbb{E}f(X_T) - \mathbb{E}f(\bar{X}_T) = hC_f + I^h(f) + o(h),$$

*where $C_f = C(U^0)$ and $I^h(f) = I^h(U^0)$ are defined as in Lemma 9 with*

$$U_s^0 = \sigma'\left(\int_{-r}^0 X_{s+u}\, d\nu(u)\right) D_s f'(X_T) + b'\left(\int_{-r}^0 X_{s+u}\, d\nu(u)\right) f'(X_T)$$

$$+ \sigma'\left(\int_{-r}^0 X_{s+u}\, d\nu(u)\right) D_s\left(\int_0^T \theta_t\, dt\right) + b'\left(\int_{-r}^0 X_{s+u}\, d\nu(u)\right) \int_s^T \theta_t\, dt$$

*and $\theta$ is the unique solution of*

$$\theta_t = \alpha^*\left(J\left(f'(X_T) + \int_0^T \theta_s\, ds\right)\right)(t) + \beta^*\left(E\left(f'(X_T) + \int_t^T \theta_s\, ds\Big|\mathcal{F}.\right)\right)(t)$$

*with*

$$\alpha^*(X)(t) = \mathbb{E}\left(\int_{\max(t-T,-r)}^0 \sigma'\left(\int_{-r}^0 X_{t-u+v}\, d\nu(v)\right) X_{t-u}\, d\nu(u)\Big|\mathcal{F}_t\right),$$

$$\beta^*(X)(t) = \mathbb{E}\left(\int_{\max(t-T,-r)}^0 b'\left(\int_{-r}^0 X_{t-u+v}\, d\nu(v)\right) X_{t-u}\, d\nu(u)\Big|\mathcal{F}_t\right).$$

Proof.   We have

$$\mathbb{E}f(X_T) - \mathbb{E}f(\bar{X}_T) = \mathbb{E}\int_0^1 f'(aX_T + (1-a)\bar{X}_T)\, da(X_T - \bar{X}_T).$$

Let $F^h = \int_0^1 f'(aX_T + (1-a)\bar{X}_T)\, da$ and $Y_t = X_t - \bar{X}_t$. We have

$$b\left(\int_{-r}^0 X_{s+u}\, d\nu(u)\right) - b\left(\int_{-r}^0 \bar{X}_{\eta(s)+\eta(u)}\, d\nu(u)\right)$$

$$= b_1^h(s)\int_{-r}^0 (X_{s+u} - \bar{X}_{\eta(s)+\eta(u)})\, d\nu(u)$$

and

$$\sigma\left(\int_{-r}^0 X_{s+u}\, d\nu(u)\right) - \sigma\left(\int_{-r}^0 \bar{X}_{\eta(s)+\eta(u)}\, d\nu(u)\right)$$

$$= \sigma_1^h(s)\int_{-r}^0 (X_{s+u} - \bar{X}_{\eta(s)+\eta(u)})\, d\nu(u),$$

where

$$(19) \quad \sigma_1^h(s) = \int_0^1 \sigma'\left(a\int_{-r}^0 X_{s+u}\, d\nu(u) + (1-a)\int_{-r}^0 \bar{X}_{\eta(s)+\eta(u)}\, d\nu(u)\right) da,$$

$$(20) \quad b_1^h(s) = \int_0^1 b'\left(a\int_{-r}^0 X_{s+u}\, d\nu(u) + (1-a)\int_{-r}^0 \bar{X}_{\eta(s)+\eta(u)}\, d\nu(u)\right) da.$$

We deduce that $Y_t$ is a solution of

$$dY_t = \alpha^h(Y)(t)\, dW_t + \beta^h(Y)(t)\, dt + dG_t^h, \qquad t \geq 0,$$

$$Y_s = 0, \qquad s \in [-r, 0],$$

with

$$(21) \quad \alpha^h(Y)(s) = \sigma_1^h(s)\int_{-r}^0 Y_{s+u}\, d\nu(u),$$

$$(22) \quad \beta^h(Y)(s) = b_1^h(s)\int_{-r}^0 Y_{s+u}\, d\nu(u),$$

$$G_t^h = \int_0^t \sigma_1^h(s)\int_{-r}^0 (\bar{X}_{s+u} - \bar{X}_{\eta(s)+\eta(u)})\, d\nu(u)\, dW_s$$

$$+ \int_0^t b_1^h(s)\int_{-r}^0 (\bar{X}_{s+u} - \bar{X}_{\eta(s)+\eta(u)})\, d\nu(u)\, ds.$$

The operators $\alpha^h$ and $\beta^h$ satisfy Assumption A1 uniformly with respect to $h$. The adjoint operators $\alpha^{h*}$ and $\beta^{h*}$ are given by

$$\alpha^{h*}(X)(s) = \mathbb{E}\left(\int_{\max(s-T,-r)}^0 \sigma_1^h(s-u)X_{s-u}\, d\nu(u)\Big|\mathcal{F}_s\right),$$

$$\beta^{h*}(X)(s) = \mathbb{E}\left(\int_{\max(s-T,-r)}^{0} b_1^h(s-u)X_{s-u}\,d\nu(u)\Big|\mathcal{F}_s\right).$$

They satisfy Assumption A2 uniformly in $h$. Note that the estimates on $D_u\alpha^{h*}$, $D_u\beta^{h*}$, and so on, follow from the boundedness of the derivatives of $b$ and $\sigma$ and the following estimates (see [6]):

$$\forall\, p\in[1,+\infty) \qquad \sup_{0\leq u,v,t\leq T}\mathbb{E}n^p(X_t,u,v)<\infty \quad \text{and}$$

$$\sup_{h>0}\sup_{0\leq u,v,t\leq T}\mathbb{E}n^p(\bar{X}_t,u,v)<\infty.$$

Using Theorem 4, we obtain

$$\mathbb{E}f(X_T) - \mathbb{E}f(\bar{X}_T) = \mathbb{E}F^hY_T = \mathbb{E}F^hG_T^h + \mathbb{E}\int_0^T\theta_s^h G_s^h\,ds,$$

with

$$\theta^h = (I-(A^h)^* - (B^h)^*)^{-1}[\alpha^{h*}(J(F^h)) + \beta^{h*}(\mathbb{E}(F^h|\mathcal{F}_.))]$$

But

$$\begin{aligned}
\mathbb{E}F^hG_T^h = {}&\mathbb{E}\int_0^T D_sF^h\sigma_1^h(s)\int_{-r}^0(\bar{X}_{s+u}-\bar{X}_{\eta(s)+\eta(u)})\,d\nu(u)\,ds\\
&+\mathbb{E}\int_0^T F^h b_1^h(s)\int_{-r}^0(\bar{X}_{s+u}-\bar{X}_{\eta(s)+\eta(u)})\,d\nu(u)\,ds
\end{aligned}$$

and

$$\begin{aligned}
\mathbb{E}\int_0^T\theta_s^h G_s^h\,ds = {}&\mathbb{E}\int_0^T D_s\left(\int_0^T\theta_t^h\,dt\right)\sigma_1^h(s)\int_{-r}^0(\bar{X}_{s+u}-\bar{X}_{\eta(s)+\eta(u)})\,d\nu(u)\,ds\\
&+\mathbb{E}\int_0^T\left(\int_s^T\theta_t^h\,dt\right)b_1^h(s)\int_{-r}^0(\bar{X}_{s+u}-\bar{X}_{\eta(s)+\eta(u)})\,d\nu(u)\,ds.
\end{aligned}$$

We end the proof applying Lemma 9 with

$$U_s^h = D_sF^h\sigma_1^h(s) + F^h b_1^h(s) + D_s\left(\int_0^T\theta_t^h\,dt\right)\sigma_1^h(s) + b_1^h(s)\int_s^T\theta_t^h\,dt.$$

Since $\lim_h\sup_{0\leq u,v,t\leq T}\mathbb{E}n^p(X_t-\bar{X}_t,u,v)=0$ for $p\geq 1$, the convergence of $U_s^h$ and its derivatives follows from Proposition 8.   $\square$

## 4. The irregular case.

4.1. *Duality for the derivatives of $Y$.* In this section we want to establish duality formulae similar to those in Theorem 4, for the derivatives of the solution of (3). This is motivated by the treatment of irregular functions of the Euler scheme, which will involve integrations by parts (see Section 4.2).

We will need regularity assumptions on the operators $\alpha$ and $\beta$. We will assume that for $Z \in \mathcal{S}_a([0, T]; H)$, $\alpha(Z)$ and $\beta(Z)$ are in $\mathbb{D}^\infty$ and that we can define recursively the operators $D^k_{s_1 \cdots s_k}\alpha$ and $D^k_{s_1 \cdots s_k}\beta$ so that, for $k > 1$ and $s_1, \ldots, s_k \in [0, T]$,

$$D_{s_1}((D^{k-1}_{s_2 \cdots s_k}\alpha)(Z)) = (D^k_{s_1 \cdots s_k}\alpha)(Z) + (D^{k-1}_{s_2 \cdots s_k}\alpha)(D_{s_1}Z),$$

$$D_{s_1}((D^{k-1}_{s_2 \cdots s_k}\beta)(Z)) = (D^k_{s_1 \cdots s_k}\beta)(Z) + (D^{k-1}_{s_2 \cdots s_k}\beta)(D_{s_1}Z).$$

Note that $D^k_{s_1 \cdots s_k}Z_t$ can be viewed as an element of $H^{d^k}$, with coordinates $D^{l_1 \cdots l_k}_{s_1 \cdots s_k}$, where the $l_i$ superscripts refer to the coordinates of $W$ with respect to which differentiation occurs ($l_i = 1, \ldots, d$). We now introduce the following assumption.

ASSUMPTION A3. For $\gamma = \alpha$ and $\gamma = \beta$, the operators $D^j_{s_1 \cdots s_j}\gamma$ defined above are bounded from $L^q_a([0, T]; H)$ into $L^p_a([0, T]; H^{d^j})$ for $2 \leq p < q < \infty$. Moreover, their adjoints satisfy the following estimates. For all positive integers $j$ and for $1 \leq p < q \leq 2$, there exists a positive constant $C_{p,q,j}$ such that for all $u, v, s_1, \ldots, s_j \in [0, T]$, we have

$$\|(D^j_{s_1 \cdots s_j}\gamma)^*\|_{q \to p} + \|D_u(D^j_{s_1 \cdots s_j}\gamma)^*\|_{q \to p} + \|D^2_{uv}(D^j_{s_1 \cdots s_j}\gamma)^*\|_{q \to p} \leq C_{p,q,j},$$

where $\|\cdot\|_{q \to p}$ stands for the operator norm from $L^q_a$ into $L^p_a$, and the operators $D_u(D^j_{s_1 \cdots s_j}\gamma)^*$, $D^2_{uv}(D^j_{s_1 \cdots s_j}\gamma)^*$ are defined in the same way as $D_u\gamma^*$ in Assumption A2.

Recall the notation

$$n(Z, u, v) = |Z| + |D_uZ| + |D_vZ| + |D^2_{uv}Z|.$$

THEOREM 11. *Let Assumptions* A1, A2 *and* A3 *hold. Let* $G = (G_t)_{0 \leq t \leq T}$ *be a continuous adapted process with values in $H$, satisfying $G_t \in \mathbb{D}^\infty$, for all $t \in [0, T]$. Then the solution of (3) satisfies $Y_t \in \mathbb{D}^\infty$, for $t \in [0, T]$.*

*Moreover, given an adapted process $F = (F_t)_{0 \leq t \leq T}$ with values in $L^2([0, T]^k; H^{d^k})$ such that $F_t \in \mathbb{D}^{2,2}$ for $t \in [0, T]$, and*

$$\sup_{0 \leq u, v, s_1, \ldots, s_k \leq T} \left( \mathbb{E} \int_0^T n^2(F_t(s_1, \ldots, s_k), u, v) \, dt \right)^{1/2} < \infty,$$

*there exist adapted processes* $\theta^{(0)}, \ldots, \theta^{(j)}, \ldots, \theta^{(k)}$, *with values in* $H, \ldots,$
$L^p([0,T]^j; H^{d^j}), \ldots, L^p([0,T]^k; H^{d^k})$ *respectively for all* $p \in [1,2)$, *such that*

$$(23) \qquad \mathbb{E} \int_0^T \langle F_t, D^k Y_t \rangle \, dt = \sum_{j=0}^k \mathbb{E} \int_0^T \langle \theta_t^{(j)}, D^j G_t \rangle \, dt,$$

*and, for* $1 \le p < q \le 2$, $j = 0, \ldots, k$,

$$(24) \qquad \begin{aligned} \sup_{u,v,s_1,\ldots,s_j \in [0,T]} &\left( \mathbb{E} \int_0^T n^p(\theta_t^{(j)}(s_1, \ldots, s_j), u, v) \, dt \right)^{1/p} \\ &\le C_{p,q,k} \sup_{u,v,s_1,\ldots,s_k \in [0,T]} \left( \mathbb{E} \int_0^T n^q(F_t(s_1, \ldots, s_k), u, v) \, dt \right)^{1/q}, \end{aligned}$$

*where the constants* $C_{p,q,k}$ *do not depend on* $F$.

COROLLARY 12. *In addition to the assumptions of Theorem* 11, *we assume that the process* $(D^k G_t)$ *is right-continuous with respect to* $t$ *for* $k \ge 1$. *Then, if* $Y$ *is the solution of* (3), *the process* $(D^k Y_t)_{0 \le t \le T}$ *has a right-continuous modification. Moreover, given a random variable* $\Phi$ *with values in* $L^2([0,T]^k; H^{d^k})$, *such that* $\Phi \in \mathbb{D}^{2,2}$ *and*

$$\sup_{u,v,s_1,\ldots,s_k \in [0,T]} \mathbb{E} n^2(\Phi(s_1, \ldots, s_k), u, v) < \infty,$$

*there exist adapted processes* $\hat{\theta}^{(0)}, \ldots, \hat{\theta}^{(j)}, \ldots, \hat{\theta}^{(k)}$, *with values in* $H, \ldots,$
$L^p([0,T]^j; H^{d^j}), \ldots, L^p([0,T]^k; H^{d^k})$ *respectively for all* $p \in [1,2)$, *such that*

$$(25) \qquad \mathbb{E} \langle \Phi, D^k Y_T \rangle = \mathbb{E} \langle \Phi, D^k G_T \rangle + \sum_{j=0}^k \mathbb{E} \int_0^T \langle \hat{\theta}_t^{(j)}, D^j G_t \rangle \, dt,$$

*and, for* $1 \le p < q \le 2$, $j = 0, \ldots, k$,

$$(26) \qquad \begin{aligned} \sup_{u,v,s_1,\ldots,s_j \in [0,T]} &\left( \mathbb{E} \int_0^T n^p(\hat{\theta}_t^{(j)}(s_1, \ldots, s_j), u, v) \, dt \right)^{1/p} \\ &\le C_{p,q,k} \sup_{u,v,s_1,\ldots,s_k \in [0,T]} (\mathbb{E} n^q(\Phi(s_1, \ldots, s_k), u, v))^{1/q}, \end{aligned}$$

*where the constants* $C_{p,q,k}$ *do not depend on* $\Phi$.

PROOF. By differentiating (3) $k$ times, we get

$$D^k Y_t = D^k G_t + \int_0^t D^k(\alpha(Y)(s)) \cdot dW_s + \int_0^t D^k(\beta(Y)(s)) \, ds + I_t^{(k)},$$

with

$$I_t^{(k)}(s_1, \ldots, s_k) = \sum_{i=1}^k D_{s_1 \cdots \hat{s_i} \cdots s_k}^{k-1}(\alpha(Y)(s_i)) \mathbb{1}_{\{s_i \le t\}},$$

where the notation $\hat{s_i}$ means that the variable $s_i$ is omitted. More precisely, recall that $D_{s_1 \cdots s_k}^k Y_t$ can be viewed as an element of $H^{d^k}$, with coordinates $D_{s_1 \cdots s_k}^{l_1 \cdots l_k}$, where the $l_i$ superscripts refer to the coordinates of $W$ with respect to which differentiation occurs. With this more precise notation, we have

$$I_t^{l_1 \cdots l_k}(s_1, \ldots, s_k) = \sum_{i=1}^k D_{s_1 \cdots \hat{s_i} \cdots s_k}^{l_1 \cdots \hat{l_i} \cdots l_k}(\alpha_{l_i}(Y)(s_i)) \mathbb{1}_{\{s_i \le t\}}.$$

We can write

$$
\begin{aligned}
(27) \qquad D_{s_1 \cdots s_k}^k Y_t = {} & G_t^{(k)}(s_1, \ldots, s_k) + \int_0^t \alpha(D_{s_1 \cdots s_k}^k Y)(s) \cdot dW_s \\
& + \int_0^t \beta(D_{s_1 \cdots s_k}^k Y)(s) \, ds,
\end{aligned}
$$

with

$$
\begin{aligned}
(28) \qquad G_t^{(k)}(s_1, \ldots, s_k) = {} & D_{s_1 \cdots s_k}^k G_t + \sum_{j=1}^k \sum_{\tau \in A_j^k} \int_0^t (D_{s_\tau}^j \alpha)(D_{s_{\bar\tau}}^{k-j} Y)(s) \cdot dW_s \\
& + \int_0^t (D_{s_\tau}^j \beta)(D_{s_{\bar\tau}}^{k-j} Y)(s) \, ds + I_t^{(k)}(s_1, \ldots, s_k),
\end{aligned}
$$

where $A_j^k$ is the set of $j$-tuples $\tau = (i_1, \ldots, i_j)$, with $1 \le i_1 < \cdots < i_j \le k$ $s_\tau = s_{i_1} \cdots s_{i_j}$, and $\bar\tau$ is the ordered complement of $\tau$. We deduce from (27) that $D_{s_1 \cdots s_k}^k Y_t$ solves an equation similar to (3) and we can derive $L^p - L^q$ estimates for the derivatives of $Y$ using Assumption A3 and the regularity of $G$. Here again the formal differentiation can be justified by a standard approximation argument.

We now prove (23) by induction. For $k = 0$, the result reduces to the first part of Theorem 4 and Proposition 6. Now assume that (23) and (24) hold up to the order $k - 1$. We first deduce from (27) that

$$\mathbb{E} \int_0^T \langle F_t, D^k Y_t \rangle \, dt = \mathbb{E} \int_0^T \langle G_t^{(k)}, \theta_t^{(k)} \rangle \, dt,$$

with $\theta_t^{(k)}(s_1, \ldots, s_k) = (I - A^* - B^*)^{-1}(F(s_1, \ldots, s_k))$. The estimates for $\theta^{(k)}$ follow from Proposition 6. It follows from (28) that

$$\mathbb{E} \int_0^T \langle G_t^{(k)}, \theta_t^{(k)} \rangle \, dt = \mathbb{E} \int_0^T \langle D^k G_t, \theta_t^{(k)} \rangle \, dt + R_T,$$

where $R_T$ is the sum of terms which are of three types, which we study successively.

TYPE 1.

$$\mathbb{E}\int_0^T dt \int_{[0,T]^k} dt_1 \cdots dt_k \left\langle \theta_t^{(k)}(t_1,\ldots,t_k), \left(\int_0^t (D_{t_1\cdots t_j}^j \alpha)(D_{t_{j+1}\cdots t_k}^{k-j} Y)(s)\cdot dW_s\right)\right\rangle,$$

with $j \geq 1$. We have

$$\mathbb{E}\int_0^T dt \int_{[0,T]^k} dt_1 \cdots dt_k \left\langle \theta_t^{(k)}(t_1,\ldots,t_k), \left(\int_0^t (D_{t_1\cdots t_j}^j \alpha)(D_{t_{j+1}\cdots t_k}^{k-j} Y)(s)\cdot dW_s\right)\right\rangle$$

$$= \int_{[0,T]^k} dt_1 \cdots dt_k$$

$$\times \mathbb{E}\int_0^T dt \left\langle \theta_t^{(k)}(t_1,\ldots,t_k), \left(\int_0^t (D_{t_1\cdots t_j}^j \alpha)(D_{t_{j+1}\cdots t_k}^{k-j} Y)(s)\cdot dW_s\right)\right\rangle$$

$$= \int_{[0,T]^k} dt_1 \cdots dt_k \, \mathbb{E}\int_0^T \langle \tilde{\theta}_s(t_1,\ldots,t_k), D_{t_{j+1}\cdots t_k}^{k-j} Y_s\rangle \, ds,$$

with

$$\tilde{\theta}_t(t_1,\ldots,t_k) = (D_{t_1\cdots t_j}^j \alpha)^* \left(J\left(\int_0^T \theta_s^{(k)}(t_1,\ldots,t_k)\,ds\right)\right)(t).$$

Here, we have used Proposition 3, with $D_{t_1\cdots t_k}^j \alpha$ instead of $\alpha$. Using Assumption A3, we have, as in the proof of Proposition 6,

$$\left(\mathbb{E}\int_0^T n^p(\tilde{\theta}_t(t_1,\ldots,t_k),u,v)\,dt\right)^{1/p}$$

$$\leq C_{p,q}\left(\mathbb{E}\int_0^T n^q(\theta_t^{(k)}(t_1,\ldots,t_k),u,v)\,dt\right)^{1/q},$$

for $1 \leq p < q \leq 2$. We can now apply the induction hypothesis, for a fixed $t_1,\ldots,t_j$, to the process $\tilde{\theta}_t(t_1,\ldots,t_j)$ considered as a process with values in $L^p([0,T]^{k-j}; H^{d^{k-j}})$. Note that these processes are not in $L^2$ but, by a suitable density argument, the induction hypothesis can be applied.

TYPE 2.

$$\mathbb{E}\int_0^T dt \int_{[0,T]^k} dt_1 \cdots dt_k \left\langle \theta_t^{(k)}(t_1,\ldots,t_k), \left(\int_0^t (D_{t_1\cdots t_j}^j \beta)(D_{t_{j+1}\cdots t_k}^{k-j} Y)(s)\,ds\right)\right\rangle,$$

with $j \geq 1$. We have

$$\mathbb{E}\int_0^T dt \int_{[0,T]^k} dt_1 \cdots dt_k \left\langle \theta_t^{(k)}(t_1,\ldots,t_k), \left(\int_0^t (D_{t_1\cdots t_j}^j \beta)(D_{t_{j+1}\cdots t_k}^{k-j} Y)(s)\,ds\right)\right\rangle$$

$$= \int_{[0,T]^k} dt_1 \cdots dt_k$$
$$\times \, \mathbb{E} \int_0^T dt \left\langle \theta_t^{(k)}(t_1, \ldots, t_k), \left( \int_0^t (D_{t_1 \cdots t_j}^j \beta)(D_{t_{j+1} \cdots t_k}^{k-j} Y)(s) \, ds \right) \right\rangle$$
$$= \int_{[0,T]^k} dt_1 \cdots dt_k \mathbb{E} \int_0^T dt \, \langle \tilde{\theta}_t(t_1, \ldots, t_k), D_{t_{j+1} \cdots t_k}^{k-j} Y_t \rangle,$$

with

$$\tilde{\theta}_t(t_1, \ldots, t_k) = (D_{t_1 \cdots t_k}^j \beta)^* \left( \mathbb{E} \left( \int_{\cdot}^T \theta_s^{(k)}(t_1, \ldots, t_k) \, ds \, \middle| \, \mathcal{F}. \right) \right)(t).$$

Here again, we have used a variant of Proposition 3, with $D_{t_1 \cdots t_k}^j \beta$ instead of $\beta$, and the $L^p - L^q$ estimate follows from Assumption A3 as in the proof of Proposition 6.

TYPE 3. The terms of type 3 come from $I^{(k)}$. They are of the following form (we let $\theta = \theta^{(k)}$):

$$\mathbb{E} \int_0^T dt \int_{[0,T]^k} ds_1 \cdots ds_k$$
$$\times \, \langle (D_{s_1 \cdots s_{j-1}}^{j-1} \alpha)(D_{s_{j+1} \cdots s_k}^{k-j} Y)(s_j) \mathbb{1}_{\{s_j \le t\}} \theta_t(s_1, \ldots, s_j, \ldots, s_k) \rangle$$

and can be treated as follows (the notation $\int_{[0,T]^{k-1}} d\hat{s_j}$ means that integration with respect to $s_j$ is omitted):

$$\mathbb{E} \int_0^T dt \int_{[0,T]^k} ds_1 \cdots ds_k$$
$$\times \, \langle (D_{s_1 \cdots s_{j-1}}^{j-1} \alpha)(D_{s_{j+1} \cdots s_k}^{k-j} Y)(s_j) \mathbb{1}_{\{s_j \le t\}}, \theta_t(s_1, \ldots, s_j, \ldots, s_k) \rangle$$
$$= \int_{[0,T]^{k-1}} d\hat{s_j}$$
$$\times \, \mathbb{E} \int_0^T ds_j \int_{s_j}^T dt \, \langle \theta_t(s_1, \ldots, s_j, \ldots, s_k), (D_{s_1 \cdots s_{j-1}}^{j-1} \alpha)(D_{s_{j+1} \cdots s_k}^{k-j} Y)(s_j) \rangle$$
$$= \int_{[0,T]^{k-1}} d\hat{s_j} \, \mathbb{E} \int_0^T ds \left\langle \mathbb{E} \left( \int_s^T \theta_u(s_1, \cdot, s, \cdot, s_k) \, du \, \middle| \, \mathcal{F}_s \right), \right.$$
$$\left. (D_{s_1 \cdots s_{j-1}}^{j-1} \alpha)(D_{s_{j+1} \cdots s_k}^{k-j} Y)(s) \right\rangle$$
$$= \int_{[0,T]^{k-1}} d\hat{s_j} \, \mathbb{E} \int_0^T \langle \tilde{\theta}_t(s_1, \ldots, s_{j-1}, s_{j+1}, \ldots, s_k), D_{s_{j+1} \cdots s_k}^{k-j} Y_t \rangle \, dt,$$

with

$$\tilde{\theta}_t(s_1, \ldots, s_{j-1}, s_{j+1}, \ldots, s_k)$$
$$= (D_{s_1 \cdots s_{j-1}}^{j-1} \alpha)^* \left( \mathbb{E} \left( \int_{\cdot}^{T} \theta_u(s_1, \ldots, s_{j-1}, \cdot, s_{j+1}, \ldots, s_k) \, du \Big| \mathcal{F}_{\cdot} \right) \right)(t).$$

We have, using Assumption A3,

$$\left( \mathbb{E} \int_0^T |\tilde{\theta}_t(s_1, \ldots, s_{j-1}, s_{j+1}, \ldots, s_k)|^p \, dt \right)^{1/p}$$

$$\leq C_{p,q} \left( \mathbb{E} \int_0^T \left| \mathbb{E} \left( \int_t^T \theta_u(s_1, \ldots, s_{j-1}, t, s_{j+1}, \ldots, s_k) \, du \Big| \mathcal{F}_t \right) \right|^q dt \right)^{1/q}$$

$$\leq C_{p,q} \left( \mathbb{E} \int_0^T \left| \int_t^T \theta_s(s_1, \ldots, s_{j-1}, t, s_{j+1}, \ldots, s_k) \, ds \right|^q dt \right)^{1/q}$$

$$\leq C_{p,q} T \sup_{0 \leq t \leq T} \left( \mathbb{E} \int_0^T |\theta_s(s_1, \ldots, s_{j-1}, t, s_{j+1}, \ldots, s_k)|^q \, ds \right)^{1/q}.$$

Finally, we have $\mathbb{E} \int_0^T \langle F_t, D^k Y_t \rangle \, dt$ as a sum of terms like $\mathbb{E} \int_0^T \langle \theta_t^{(j)}, D^j G_t \rangle \, dt$, with $j \leq k$, or $\mathbb{E} \int_0^T \langle F_t^{(i)}, D^i Y_t \rangle \, dt$, with $i \leq k-1$, with appropriate estimates for the processes $F^{(i)}$. We can now apply the induction hypothesis to terms like $\mathbb{E} \int_0^T \langle F_t^{(i)}, D^i Y_t \rangle \, dt$. The various $\theta$'s given by the induction hypothesis combine to produce the final form stated in Theorem 11.  □

PROOF OF COROLLARY 12.  The continuity of the derivatives of $Y$ follow easily from the assumptions on $G$ and (27). Using the notation of the proof of Theorem 11, we have

$$\mathbb{E} \langle \Phi, D^k Y_T \rangle = \mathbb{E} \langle \Phi, G_T^{(k)} \rangle + \mathbb{E} \left\langle \Phi, \int_0^T \alpha(D^k Y)(s) \cdot dW_s + \int_0^T \beta(D^k Y)(s) \, ds \right\rangle$$

$$= \mathbb{E} \langle \Phi, G_T^{(k)} \rangle + \mathbb{E} \int_0^T \langle \alpha^*(J(\Phi))(s), D^k Y_s \rangle \, ds$$

$$+ \mathbb{E} \int_0^T \langle \beta^*(\mathbb{E}(\Phi|\mathcal{F}_{\cdot}))(s), D^k Y_s \rangle \, ds.$$

For the last two terms, we can apply Theorem 11, and the fact that, for $1 \leq p < q \leq 2$,

$$\left( \mathbb{E} \int_0^T n^p(J_t(\Phi(s_1 \cdots s_k)), u, v) \, dt \right)^{1/p} \leq C_{p,q} (\mathbb{E} n^q(\Phi(s_1 \cdots s_k), u, v))^{1/q}$$

and

$$\left( \mathbb{E} \int_0^T n^p(\mathbb{E}(\Phi(s_1 \cdots s_k)|\mathcal{F}_t), u, v) \, dt \right)^{1/p} \leq C_{p,q} (\mathbb{E} n^q(\Phi(s_1 \cdots s_k), u, v))^{1/q}.$$

For the first term, we have, using (28),

$$\mathbb{E}\langle \Phi, G_T^{(k)}\rangle = \mathbb{E}\langle \Phi, D^k Y_T\rangle + R_T,$$

where $R_T$ is a sum of terms of three different types, which can be treated in the same way as in the proof of Theorem 11. □

4.2. *Application to the Euler approximation of delay equations*: *the irregular case*. In this section we consider the processes $X$ and $\bar{X}$ of Section 3.4. We assume that $f$ is a measurable bounded function, $b$ and $\sigma$ are in $\mathcal{C}_b^\infty$ and that the variable $X_T$ is nondegenerate:

$$(29) \qquad \mathbb{E}((\gamma_{X_T})^{-p}) < \infty \qquad \forall p > 1,$$

where $\gamma_{X_T}$ denotes the Malliavin covariance matrix of $X_T$. This condition is satisfied in the uniformly elliptic case, that is, $\sigma(x) \geq a > 0$ for all $x \in \mathbb{R}$ (see [10]) and under weaker assumptions (see [2]) when $\nu$ is a Dirac measure.

THEOREM 13. *For $b$ and $\sigma$ in $\mathcal{C}_b^\infty$ and $X_T$ satisfying* (29), *we have for $f$ measurable bounded*:

$$|\mathbb{E}f(X_T) - \mathbb{E}f(\bar{X}_T)| \leq C_f h.$$

SKETCH OF PROOF. We consider the following truncation function. Let $\Psi : [0, +\infty) \mapsto \mathbb{R}$ be a $\mathcal{C}^\infty$ function with bounded derivatives such that $\mathbb{1}_{[0,1/8]} \leq \Psi \leq \mathbb{1}_{[0,1/4]}$ and let $\gamma_{X_T}$ be the Malliavin covariance matrix of $X_T$ be which in our one-dimensional setting reduces to $\gamma_{X_T} = \int_0^T (D_u X_T)^2 \, du$. We define $\Psi_T^h$ by

$$(30) \qquad \Psi_T^h = \Psi\left(\frac{\int_0^T (D_u X_T - D_u \bar{X}_T)^2 \, du}{\gamma_{X_T}}\right).$$

Observe that

$$\mathbb{P}(\Psi_T^h \neq 1) \leq \mathbb{P}\left(\gamma_{X_T}^{-1} \int_0^T (D_u X_T - D_u \bar{X}_T)^2 \, du > 1/8\right)$$

$$\leq 8^p \mathbb{E}\gamma_{X_T}^{-p}\left(\int_0^T (D_u X_T - D_u \bar{X}_T)^2 \, du\right)^p,$$

for all $p > 0$. Using Hölder's inequality and $\sup_u \mathbb{E}|D_u X_T - D_u \bar{X}_T|^p \leq Ch^{p/2}$ one can easily prove that for all $p \geq 1$

$$(31) \qquad \mathbb{P}(\Psi_T^h \neq 1) \leq Ch^p.$$

Moreover, we have

$$\{\Psi_T^h \neq 0\} \subset \left\{\int_0^T (D_u X_T - D_u \bar{X}_T)^2 \, du \leq \gamma_{X_T}/4\right\}.$$

Since

$$\gamma_{aX_T+(1-a)\bar{X}_T} \geq \gamma_{X_T}/2 - (1-a)^2 \int_0^T (D_u X_T - D_u \bar{X}_T)^2 \, du,$$

we obtain for $0 \leq a \leq 1$

(32)                    $$\{\Psi_T^h \neq 0\} \subset \{\gamma_{aX_T+(1-a)\bar{X}_T} \geq \tfrac{1}{4}\gamma_{X_T}\}.$$

Now we have

$$\mathbb{E}f(X_T) - \mathbb{E}f(\bar{X}_T)$$
$$= \mathbb{E}\{(f(X_T) - f(\bar{X}_T))(1 - \Psi_T^h)\} + \mathbb{E}\{(f(X_T) - f(\bar{X}_T))\Psi_T^h\}.$$

By construction of $\Psi_T^h$, the first term is of order $h^p$ for all $p \geq 1$ and we just have to prove that the second one is of order $h$.

Let $Y_T = X_T - \bar{X}_T$ and let $(f_m)$ be a sequence of $\mathcal{C}^1$ functions such that $\|f_m\|_\infty \leq \|f\|_\infty$ and $(f_m)$ converges $dx$ a.e. to $f$. We have

$$\mathbb{E}(f_m(X_T) - f_m(\bar{X}_T))\Psi_T^h = \int_0^1 \mathbb{E}f_m'(aX_T + (1-a)\bar{X}_T)Y_T\Psi_T^h \, da.$$

Since $X_T$ admits a density $\mathbb{E}f_m(X_T)\Psi_T^h$ converges to $\mathbb{E}f(X_T)\Psi_T^h$. Now on the set $\{\Psi_T^h \neq 0\}$, $\det \gamma_{\bar{X}_T} > 0$ and from [13], Corollary 2.2.1, page 88 (see [3] in higher dimension) $\bar{X}_T$ has an absolutely continuous law conditioned by $\{\Psi_T^h \neq 0\}$ and $\mathbb{E}f_m(\bar{X}_T)\Psi_T^h$ converges to $\mathbb{E}f(\bar{X}_T)\Psi_T^h$. It remains to prove that $\mathbb{E}(f_m(X_T) - f_m(\bar{X}_T))\Psi_T^h$ is bounded by $C_f h$ where the constant $C_f$ only depends on $f$ through $\|f\|_\infty$. Using the Malliavin integration by parts formula (see [14]) we obtain

$$\mathbb{E}(f_m(X_T) - f_m(\bar{X}_T))\Psi_T^h$$
$$= \int_0^1 \mathbb{E}g_m(aX_T + (1-a)\bar{X}_T)H_3(aX_T + (1-a)\bar{X}_T, Y_T\Psi_T^h) \, da,$$

where $g_m(x) = \int_0^x dy \int_0^y dz \, f_m(z)$ so that $g_m$ is in $\mathcal{C}^3$ and $g_m'' = f_m$ and $H_3$ is defined recursively by

$$H_1(F, G) = G\gamma_F^{-1}\delta(DF) - \langle D(G\gamma_F^{-1}), DF \rangle,$$
$$H_k(F, G) = H_1(F, H_{k-1}(F, G)), \qquad k \geq 2,$$

with

$$\langle DH, DF \rangle = \int_0^T D_u H D_u F \, du.$$

Hence for any measurable set $A$, we have (see [1])

(33)                $$\|H_k(F, G)\mathbb{1}_A\|_p \leq C\|\gamma_F^{-1}\mathbb{1}_A\|_{p_1}^{k_1}\|F\|_{p_2, q_2}^{k_2}\|G\|_{p_3, q_3}$$

for some constants $C$, $k_1$, $p_1$, $k_2$, $p_2$, $q_2$, $p_3$, $q_3$, depending on $k$ and $p$.

Observe that

$$H_3(aX_T + (1-a)\bar{X}_T, Y_T \Psi_T^h) = \sum_{i=0}^{3} \langle \Phi_i^h, D^i Y_T \rangle,$$

for smooth variables $\Phi_i^h$, and finally

$$\mathbb{E}(f_m(X_T) - f_m(\bar{X}_T))\Psi_T^h = \mathbb{E}\sum_{i=0}^{3} \left\langle \int_0^1 g_m(aX_T + (1-a)\bar{X}_T)\Phi_i^h \, da, D^i Y_T \right\rangle$$

$$= \mathbb{E}\sum_{i=0}^{3} \langle F_m^{i,h}, D^i Y_T \rangle.$$

Moreover, the variables $F_m^{i,h} \in \mathbb{D}^{2,2}$ and from (33) and (39) we deduce the following estimate:

$$\sup_{u,v} \mathbb{E}n^2(F_m^{i,h}, u, v) \le C\|\gamma_{X_T}^{-1}\|_{p_1}^{k_1}\|aX_T + (1-a)\bar{X}_T\|_{p_2,q_2}^{k_2},$$

for some $k_1$, $p_1$, $k_2$, $p_2$, $q_2$, with $C$ independent of $m$. Applying Corollary 12, this gives

$$\mathbb{E}(f_m(X_T) - f_m(\bar{X}_T))\Psi_T^h = \sum_{i=0}^{3} \left( \mathbb{E}\langle F_m^{i,h}, D^i G_T \rangle + \sum_{j=0}^{i} \mathbb{E}\int_0^T \langle \hat{\theta}_t^{(i,j)}, D^i G_t \rangle \, dt \right),$$

where the process $G_t^h$ is defined in Section 3.4 by

$$G_t^h = \int_0^t \sigma_1^h(s) \int_{-r}^0 (\bar{X}_{s+u} - \bar{X}_{\eta(s)+\eta(u)}) \, d\nu(u) \, dW_s$$

$$+ \int_0^t b_1^h(s) \int_{-r}^0 (\bar{X}_{s+u} - \bar{X}_{\eta(s)+\eta(u)}) \, d\nu(u) \, ds$$

$$= \int_0^t \sigma_1^h(s) \int_{-r}^0 \int_{\eta(s)+\eta(u)}^{s+u} \tilde{\sigma}_{\eta(v)}^h \, dW_v \, d\nu(u) \, dW_s$$

$$+ \int_0^t b_1^h(s) \int_{-r}^0 \int_{\eta(s)+\eta(u)}^{s+u} \tilde{b}_{\eta(v)}^h \, dv \, d\nu(u) \, ds.$$

$\sigma_1^h$, $b_1^h$, $\tilde{\sigma}_v^h$ and $\tilde{b}_v^h$ are, respectively, defined in (19), (20), (14) and (15). We end the proof using the duality relationship as in Lemma 9.  $\square$

REMARK 14. The above proof carries over to a multidimensional setting. The only technical difficulty is to extend the localization argument. Suppose $X_T$ takes its values in $\mathbb{R}^r$ with $r > 1$, and the Malliavin covariance matrix of $X_T$ is defined by $(\gamma_{X_T})_{i,j} = \langle DX_T^i, DX_T^j \rangle$, $1 \le i \le j \le r$. We want to define a smooth functional $\Psi_T^h$ such that outside a set where all the Malliavin derivatives $D^k \Psi_T^h$ vanish, we have a uniform control of the determinant of the Malliavin matrix of $aX_T + (1-a)\bar{X}_T$, $0 \leq a \leq 1$. Let $\Psi : [0, +\infty) \to \mathbb{R}$ be a $\mathcal{C}^\infty$ function with bounded derivatives such that $\mathbb{1}_{[0,1/8]} \leq \Psi \leq \mathbb{1}_{[0,1/4]}$. We define $\Psi_T^h$ by

$$(34) \qquad \Psi_T^h = \Psi\left( \frac{|D(X_T - \bar{X}_T)|^2 (1 + \|\gamma_{X_T}\|_2^2)^{(r-1)/2}}{\det \gamma_{X_T}} \right),$$

where $|DX_T|^2 = \sum_i |DX_T^i|^2 = \sum_{i,k} \int_0^T (D_u^k X_T^i)^2 \, du$ and $\|\gamma_{X_T}\|_2$ is the Hilbert–Schmidt norm of the matrix $\gamma_{X_T}$, that is, the $l^2$ norm of the coefficients. We denote by $\|\gamma_{X_T}\|$ the operator norm.

Note that $\Psi_T^h \in \mathbb{D}^\infty$ [the sum of squares of the coefficients $(\gamma_{X_T})_{i,j}$ is smooth] and that

$$\{\Psi_T^h \neq 1\} \subset \left\{ |D(X_T - \bar{X}_T)|^2 \geq \frac{\det \gamma_{X_T}}{8(1 + \|\gamma_{X_T}\|_2^2)^{(r-1)/2}} \right\},$$

so that for all $p \geq 1$, we have as in the one-dimensional case

$$(35) \qquad \exists\, C > 0 \qquad \mathbb{P}(\Psi_T^h \neq 1) \leq Ch^p,$$

provided $\|D(X_T - \bar{X}_T)\|_p \leq C_p \sqrt{h}$ for all $p \in [1, +\infty)$.

Now observe that we have the following inequality for any positive-definite $r$-dimensional matrix $A$:

$$(36) \qquad \|A\| \leq \|A\|_2 \leq \sqrt{r}\|A\|.$$

Moreover, if $\lambda_1(A)$ is the smallest eigenvalue of $A$, we have

$$(37) \qquad \lambda_1(A)^r \leq \det A \leq \lambda_1(A)\|A\|^{r-1}.$$

Observe that $\lambda_1(A) = \inf_{|\xi|=1} \xi^t A \xi$, where $|\xi|$ is the Euclidean norm of $\xi$ in $\mathbb{R}^r$. Now for $a \in [0,1]$, we derive a uniform lower bound for the smallest eigenvalue of $\gamma_{X_T + (1-a)\bar{X}_T}$:

$$
\begin{aligned}
(38) \qquad \sqrt{\lambda_1(\gamma_{X_T + (a-1)(X_T - \bar{X}_T)})} &= \inf_{|\xi|=1} \sqrt{\xi^t \gamma_{X_T + (a-1)(X_T - \bar{X}_T)} \xi} \\
&= \inf_{|\xi|=1} \left| \sum_i \xi_i D(X_T^i + (a-1)(X_T^i - \bar{X}_T^i)) \right| \\
&\geq \inf_{|\xi|=1} \sqrt{\xi^t \gamma_{X_T} \xi} - \sup_{|\xi|=1} \left| \sum_i \xi_i D(X_T^i - \bar{X}_T^i) \right| \\
&\geq \sqrt{\lambda_1(\gamma_{X_T})} - |D(X_T - \bar{X}_T)|.
\end{aligned}
$$

Moreover, we have

$$\bigcup_k \{D^k \Psi_T^h \neq 0\} \subset \left\{ |D(X_T - \bar{X}_T)|^2 \leq \frac{\det \gamma_{X_T}}{4(1 + \|\gamma_{X_T}\|_2^2)^{(r-1)/2}} \right\}.$$

But from (36) and (37) we have

$$\frac{\det \gamma_{X_T}}{4(1 + \|\gamma_{X_T}\|_2^2)^{(r-1)/2}} \leq \frac{\lambda_1(\gamma_{X_T})}{4}.$$

We deduce that

$$\bigcup_k \{D^k \Psi_T^h \neq 0\} \subset \left\{ |D(X_T - \bar{X}_T)|^2 \leq \frac{\lambda_1(\gamma_{X_T})}{4} \right\}.$$

Now if $|D(X_T - \bar{X}_T)|^2 \leq \lambda_1(\gamma_{X_T})/4$, it follows from (39) that

$$\sqrt{\lambda_1(\gamma_{X_T + (1-a)(X_T - \bar{X}_T)})} \geq \sqrt{\lambda_1(\gamma_{X_T})}/2$$

and consequently, using (37),

$$\det \gamma_{X_T + (1-a)(X_T - \bar{X}_T)} \geq \lambda_1(\gamma_{X_T})^r / 4^r \geq (\det \gamma_{X_T})^r / (\|\gamma_{X_T}\|^{r(r-1)} 4^r).$$

Finally we obtain for $0 \leq a \leq 1$

$$(39) \qquad \bigcup_k \{D^k \Psi_T^h \neq 0\} \subset \left\{ \det \gamma_{aX_T + (1-a)\bar{X}_T} \geq \frac{(\det \gamma_{X_T})^r}{\|\gamma_{X_T}\|^{r(r-1)} 4^r} \right\}.$$

Therefore on the set $\bigcup_k \{D^k \Psi_T^h \neq 0\}$, we have a control of the determinant of the inverse of the Malliavin matrix of $aX_T + (1-a)\bar{X}_T$ by the random variable $\frac{\|\gamma_{X_T}\|^{r(r-1)}}{(\det \gamma_{X_T})^r}$.

**Aknowledgments.** This research originated from discussions with Vlad Bally. The authors are grateful for his special contribution and encouragement. A. Kohatsu-Higa thanks the support of the Mathfi project (INRIA) and the Department of Mathematics of the University Marne-la-Vallée for their hospitality during various visits.

## REFERENCES

[1] BALLY, V. and TALAY, D. (1996). The law of the Euler scheme for stochastic differential equations. I. Convergence rate of the distribution function. *Probab. Theory Related Fields* **104** 43–60. MR1367666

[2] BELL, D. R. and MOHAMMED, S. E. A. (1991). The Malliavin calculus and stochastic delay equations. *J. Funct. Anal.* **99** 75–99. MR1120914

[3] BOULEAU, N. and HIRSCH, F. (1991). *Dirichlet Forms and Analysis on Wiener Space.* de Gruyter, Berlin. MR1133391

[4] BUCKWAR, E. and SHARDLOW, T. (2005). Weak approximation of stochastic differential delay equations. *IMA J. Numer. Anal.* **25** 57–86. MR2110235

[5] GOBET, E. and MUNOS, R. (2005). Sensitivity analysis using Itô–Malliavin calculus and martingales, and application to stochastic optimal control. *SIAM J. Control Optim.* **43** 1676–1713. MR2137498

[6] HU, Y., MOHAMMED, S. and YAN, F. (2001). Discrete time approximations of stochastic differential systems with memory. Preprint.

[7] KOHATSU-HIGA, A. and PETTERSSON, R. (2002). Variance reduction methods for simulation of densities on Wiener space. *SIAM J. Numer. Anal.* **40** 431–450. MR1921664

[8] KOHATSU-HIGA, A. and PROTTER, P. (1994). The Euler scheme for SDE's driven by semimartingales. In *Stochastic Analysis on Infinite-Dimensional Spaces* (H. Kunita and H.-H. Kuo, eds.) 141–151. Longman Sci. Tech., Harlow. MR1415665

[9] KÜCHLER, U. and PLATEN, E. (2002). Weak discrete time approximation of stochastic differential equations with time delay. *Math. Comput. Simulation* **59** 497–507. MR1917820

[10] KUSUOKA, S. and STROOCK, D. (1982). Applications of Malliavin calculus. I. In *Stochastic Analysis (Katata/Kyoto, 1982)* 271–306. North-Holland, Amsterdam. MR0780762

[11] LAPEYRE, B. and TEMAM, E. (2001). Competitive Monte Carlo methods for the pricing of Asian options. *J. Comput. Finance* **5** 39–59.

[12] MOHAMMED, S.-E. A. (1998). Stochastic differential systems with memory: Theory, examples and applications. In *Stochastic Analysis and Related Topics VI* (L. Decreusefond, J. Gjerde, B. Øksendal and A. S. Ustunel, eds.) 1–77. Birkhäuser, Boston. MR1652338

[13] NUALART, D. (1995). *The Malliavin Calculus and Related Topics.* Springer, Berlin. MR1344217

[14] NUALART, D. (1998). Analysis on Wiener space and anticipating stochastic calculus. *Lectures on Probability Theory and Statistics. Lecture Notes in Math.* **1690** 123–227. Springer, Berlin. MR1668111

[15] TALAY, D. and TUBARO, L. (1990). Expansion of the global error for numerical schemes solving stochastic differential equations. *Stochastic Anal. Appl.* **8** 483–509. MR1091544

E. CLÉMENT
D. LAMBERTON
LABORATOIRE D'ANALYSE ET DE
MATHÉMATIQUES APPLIQUÉES
ET PROJET MATHFI
UNIVERSITÉ DE MARNE-LA-VALLÉE
5 BLD DESCARTES, CHAMPS-SUR-MARNE
77454 MARNE-LA-VALLÉE CEDEX 2
FRANCE
E-MAIL: emmanuelle.clement@univ-mlv.fr

A. KOHATSU-HIGA
GRADUATE SCHOOL
   OF ENGINEERING SCIENCES
OSAKA UNIVERSITY
MACHIKANEYAMA CHO 1-3
OSAKA 560-8531
JAPAN